\documentclass[aop, preprint ,linenumbers, 10pt]{imsart}

\setattribute{journal}{name}{}

\usepackage[OT1]{fontenc}	
\usepackage{amsmath, amsthm, amssymb}

\usepackage{upgreek}
\usepackage{mathrsfs}

\usepackage{pdfsync}
\usepackage[a4paper,top=3cm,bottom=3cm,left=2.5cm,right=2.5cm, bindingoffset=0mm]{geometry}

\usepackage[usenames,dvipsnames]{xcolor}
\usepackage{enumerate}								
\usepackage[sort,numbers]{natbib}								
\usepackage[all]{xy}									
\usepackage{verbatim}								
\usepackage[font={footnotesize,it}]{caption}
\usepackage{dsfont}									

\usepackage{multirow}

\usepackage[pdftex,bookmarks,
						hidelinks,					
						breaklinks,unicode,
						citecolor=OliveGreen,
						linkcolor=Maroon
						]{hyperref} 

\usepackage{ifthen}

\usepackage{mathtools}

\usepackage{caption}
\usepackage{subcaption}

\usepackage{tikz}
\usepackage{tikz-cd}
\usetikzlibrary{calc, positioning, arrows, patterns}
\tikzcdset{arrow style=tikz}
\tikzset{>=latex}

\definecolor{grey1}{RGB}{71,71,71}

\usepackage[smalltableaux]{ytableau}

\arxiv{}

\startlocaldefs

\newcommand{\mfa}{{\mathfrak a}}

\newcommand{\mff}{{\mathfrak f}}
\newcommand{\mfg}{{\mathfrak g}}
\newcommand{\mfh}{{\mathfrak h}}

\newcommand{\mfl}{{\mathfrak l}}

\newcommand{\mfn}{{\mathfrak n}}

\newcommand{\mfr}{{\mathfrak r}}
\newcommand{\mfs}{{\mathfrak s}}

\newcommand{\mfN}{{\mathfrak N}}

\newcommand{\mfP}{{\mathfrak P}}

\newcommand{\mfS}{{\mathfrak S}}

\newcommand{\mfU}{{\mathfrak U}}

\newcommand{\mcB}{{\mathcal B}}
\newcommand{\mcC}{{\mathcal C}}
\newcommand{\mcD}{{\mathcal D}}

\newcommand{\mcG}{{\mathcal G}}

\newcommand{\mcL}{{\mathcal L}}

\newcommand{\mcO}{{\mathcal O}}
\newcommand{\mcP}{{\mathcal P}}

\newcommand{\mcS}{{\mathcal S}}

\newcommand{\msF}{{\mathscr F}}
\newcommand{\msG}{{\mathscr G}}

\newcommand{\msP}{{\mathscr P}}

\newcommand{\mbfE}{{\mathbf E}}

\newcommand{\mbfG}{{\mathbf G}}

\newcommand{\mbfL}{{\mathbf L}}

\newcommand{\mbfP}{{\mathbf P}}

\newcommand{\mbfW}{{\mathbf W}}
\newcommand{\mbfX}{{\mathbf X}}
\newcommand{\mbfY}{{\mathbf Y}}

\newcommand{\mbfb}{{\mathbf b}}

\newcommand{\mbfe}{{\mathbf e}}

\newcommand{\mbfk}{{\mathbf k}}

\newcommand{\mbfm}{{\mathbf m}}
\newcommand{\mbfn}{{\mathbf n}}

\newcommand{\mbfp}{{\mathbf p}}

\newcommand{\mbfs}{{\mathbf s}}
\newcommand{\mbft}{{\mathbf t}}
\newcommand{\mbfu}{{\mathbf u}}

\newcommand{\mbfy}{{\mathbf y}}
\newcommand{\mbfz}{{\mathbf z}}

\newcommand{\boldalpha}{{\boldsymbol \alpha}}
\newcommand{\boldbeta}{{\boldsymbol \beta}}

\newcommand{\boldlambda}{{\boldsymbol \lambda}}

\newcommand{\mbbE}{{\mathbb E}}

\newcommand{\mbbP}{{\mathbb P}}
\newcommand{\mbbS}{{\mathbb S}}

\newcommand{\eps}{\varepsilon}

\newcommand{\defeq}{\eqqcolon}

\newcommand{\acts}{\raisebox{\depth}{\scalebox{1}[-1]{$\circlearrowleft$}}\,}

\newcommand{\emparg}{{\,\cdot\,}}

\DeclareMathOperator{\Diffeo}{Diff}
\DeclareMathOperator{\Iso}{Iso}

\newcommand{\Dir}[1]{{D}_{#1}}								
\newcommand{\IDir}[1]{{\tilde D}_{#1}}						

\DeclareMathOperator{\ev}{ev}
\DeclareMathOperator{\bd}{bd}

\DeclareMathOperator{\pr}{pr}

\newcommand{\forallae}[1]{{\textrm{\,for ${#1}$-a.e.\,}}}

\newcommand{\FT}[1]{\widehat{#1}}

\DeclareMathOperator{\Id}{Id}
\DeclareMathOperator{\id}{id}

\DeclareMathOperator{\diag}{diag}

\DeclareMathOperator{\End}{End}

\DeclareMathOperator{\eqdef}{\coloneqq}

\let\epsilon\varepsilon
\let\subset\subseteq
\let\supset\supseteq

\newcommand{\rar}{\rightarrow}

\newcommand{\klim}{\lim_{k }}
\newcommand{\hlim}{\lim_{h }}

\DeclareMathOperator{\supp}{supp}								
\DeclareMathOperator{\sgn}{sgn}								
\newcommand{\diff}{\mathop{}\!\mathrm{d}}						
	
\newcommand{\abs}[1]{\left\lvert#1\right\rvert}						
\newcommand{\norm}[1]{\left\lVert#1\right\rVert}					
\newcommand{\set}[1]{\left\{#1\right\}}							
\newcommand{\tset}[1]{\{#1\}}							
\newcommand{\ceiling}[1]{\left\lceil#1\right\rceil}					
\newcommand{\floor}[1]{\left\lfloor#1\right\rfloor}					
\newcommand{\tonde}[1]{\left(#1\right)}							
\newcommand{\quadre}[1]{\left[#1\right]}							
\newcommand{\scalar}[2]{\left\langle #1 \,\middle |\, #2\right\rangle}		

\DeclareMathOperator{\interior}{int}								

\DeclareMathOperator{\cl}{cl}									
\DeclareMathOperator{\diam}{diam}								
\newcommand{\seq}[1]{\tonde{#1}}								
\newcommand{\tseq}[1]{(#1)}								

\DeclareMathOperator{\Meas}{\mathscr M}

\newcommand{\Mb}{\mathscr M_b}
\newcommand{\Mbp}{\mathscr M_b^+}

\DeclareMathOperator{\Prob}{\mathbb P}

\newcommand{\pfwd}{\sharp}

\DeclareMathOperator{\car}{\mathds 1}

\DeclareMathOperator{\emp}{\varnothing} 

\DeclareMathOperator{\N}{{\mathbb N}}
\DeclareMathOperator{\R}{{\mathbb R}}

\renewcommand{\C}{{\mathbb C}} 	
\DeclareMathOperator{\Z}{{\mathbb Z}}

\DeclareMathOperator{\im}{im}
\newcommand{\restr}{\big\lvert}

\newcommand{\functionnn}[5]{\begin{align*}#1\colon#2&\longrightarrow#3\\ #4&\longmapsto#5\end{align*}}		

\usepackage{booktabs}

\allowdisplaybreaks

\usetikzlibrary{shapes.misc}
\tikzset{cross/.style={cross out, draw=black, minimum size=2*(#1-\pgflinewidth), inner sep=0pt, outer sep=0pt},
cross/.default={4pt}}

\usepackage{youngtab}

\newcommand\catalannumber[4]{
  (#1)
  \foreach \dir in {#4}{
    \ifnum\dir=0
    -- ++(1,0)
    \else
    -- ++(0,1)
    \fi
  } |- (#1);
  \draw[help lines] (#1) grid +(#2,#3);
  \coordinate (prev) at (#1);
  \foreach \dir in {#4}{
    \ifnum\dir=0
    \coordinate (dep) at (1,0);
    \else
    \coordinate (dep) at (0,1);
    \fi
    \draw[line width=1pt] (prev) -- ++(dep) coordinate (prev);
  };

}

\newcommand{\iref}[1]{($\ref{#1}$)}

\newcommand{\comm}{\,\,\textrm{,}\quad\,}
\newcommand{\semicolon}{\,\,\textrm{;}\quad\,}
\newcommand{\fstop}{\,\,\textrm{.}}

\newcommand{\Mnom}[1]{M_2\!\tonde{#1}}

\DeclareMathOperator{\zero}{{\mathbf 0}}
\DeclareMathOperator{\uno}{{\mathbf 1}}
\newcommand{\imu}{\mathrm{i}}

\newcommand{\av}[1]{\left\langle#1\right\rangle}

\usepackage{scrextend}						

\DeclareMathOperator{\EGF}{\mbfG_{\exp}}
\DeclareMathOperator{\GF}{\mbfG}

\newcommand{\prid}{\mathrel{\ooalign{$\lneq$\cr\raise.22ex\hbox{$\lhd$}\cr}}}

\newcommand{\compo}{\diamond}
\newcommand{\contra}{\star}
\newcommand{\contraplus}{+}

\newcommand{\Poch}[2]{\left\langle{#1}\right\rangle_{#2}}

\usepackage{stmaryrd}

\newcommand{\length}[1]{{#1}_\bullet}

\newcommand{\bag}[1]{\llbracket{#1}\rrbracket}

\renewcommand{\P}{{\mathbb P}}

\let\temp\phi
\let\phi\varphi
\let\varphi\temp

\newcommand{\Alpha}{\mathrm{A}}
\newcommand{\Beta}{\mathrm{B}}

\newcommand{\Giry}{{\msG}}
\newcommand{\DF}{{\mcD}}
\newcommand{\GP}{{\mcG}}
\newcommand{\PP}{{\mcP}}

\newcommand{\mfsl}{{\mfs\mfl}}

\numberwithin{equation}{section}
\theoremstyle{plain}
\newtheorem{thm}{Theorem}[section]
\newtheorem*{thm*}{Theorem}

\newtheorem{prop}[thm]{Proposition}
\newtheorem{lem}[thm]{Lemma}
\newtheorem{cor}[thm]{Corollary}

\theoremstyle{definition}
\newtheorem{defs}[thm]{Definition}

\theoremstyle{remark}
\newtheorem{rem}[thm]{Remark}
\endlocaldefs

\begin{document}

\begin{frontmatter}
\title{Characteristic functionals of Dirichlet measures\thanksref{T1}
}
\runtitle{Dirichlet characteristic functionals}

\begin{aug}
\author{\fnms{Lorenzo} \snm{Dello Schiavo}\thanksref{t3}\ead[label=e1]{delloschiavo@iam.uni-bonn.de}}

\thankstext{T1}{Research supported by the Collaborative Research Center 1060 and the Hausdorff Center for Mathematics.}

\thankstext{t3}{I~gratefully acknowledge stimulating discussions with Prof.s \mbox{K.-T.}~Sturm,~E.~W. Lytvynov, C.~Stroppel and M.~Gordina. I~am also indebted to Prof.s W.~Miller, A.~M.~Vershik and M. Piccioni for respectively pointing out the references~\cite{Mil73},~\cite{Ver07} and~\cite{LetPic18}.}
\runauthor{L. Dello Schiavo}

\affiliation{Universit\"at Bonn
}

\address{Institut f\"ur Angewandte Mathematik\\
Rheinische Friedrich-Wilhelms-Universit\"at Bonn\\
Endenicher Allee 60\\
DE 53115 Bonn\\
Germany\\
\printead{e1}
}

\end{aug}

\begin{abstract}
We compute characteristic functionals of Dirichlet--Ferguson measures over a locally compact Polish space and prove continuous dependence of the random measure on the parameter measure.
In finite dimension, we identify the dynamical symmetry algebra of the characteristic functional of the Dirichlet distribution with a simple Lie algebra of type~$A$. We study the lattice determined by characteristic functionals of categorical Dirichlet posteriors, showing that it has a natural structure of weight Lie algebra module and providing a probabilistic interpretation. A partial generalization to the case of the Dirichlet--Ferguson measure is also obtained.
\end{abstract}

\vspace{.5cm}
\today
\vspace{.5cm}



\end{frontmatter}
%
%
%
%
%
%

\section{Introduction and main results}\label{s:Intro}
Let~$X$ be a locally compact Polish space with Borel $\sigma$-algebra~$\mcB(X)$ and let~$\msP(X)$ be the space of probability measures on~$(X,\mcB(X))$. For~$\sigma\in \msP(X)$ we denote by~$\DF_\sigma$ the \emph{Dirichlet--Ferguson measure}~\cite{Fer73} on~$\msP(X)$ with {probability} intensity~$\sigma$.

The characteristic functional of~$\DF_\sigma$ is commonly recognized as hardly tractable~\cite{JiaDicKuo04} and any approach to~$\DF_\sigma$ based on characteristic functional methods appears de facto ruled out in the literature. Notably, this led to the introduction of different characterizing transforms (e.g. the {Markov--Krein transform}~\cite{KerTsi01,TsiVerYor00} or the {$c$-transform}~\cite{JiaDicKuo04}), inversion formulas based on characteristic functionals of other random measures (in particular, the Gamma measure, as in~\cite{RegGugDiN02}),  and, at least in the case~\mbox{$X=\R$}, to the celebrated {Markov--Krein identity}
(see e.g.~\cite
{LijReg04}).

These investigations are based on complex analysis techniques and integral representations of special functions, in particular the Lauricella hypergeometric function~${}_kF_D$~\cite{Lau1893} 
and Carlson's~$R$ function~\cite{Car77}.
The novelty in this work consists in the combinatorial/algebraic approach adopted, allowing for broader generality and far reaching connections, especially with Lie algebra theory.

\paragraph{Fourier analysis} Denote by~$\Dir{\boldalpha_k}$ the \emph{Dirichlet distribution} on the standard simplex~$\Delta^{k-1}$ with parameter~$\boldalpha_k\in \R_+^k$, which we regard as the discretization of~$\DF_\sigma$ induced by a measurable $k$-partition~$\mbfX_k$ of~$X$ (see \S\ref{s:Prelim} below).
Our first result is the following.
\begin{thm}[see Thm.~\ref{t:UnoBis}]\label{t:Uno}
The characteristic functional~$\FT{\DF_\sigma}$ of~$\DF_\sigma$ is ---~for suitable sequences of partitions~$\mbfX_k$~--- the limit of the discrete $\DF_\sigma$-martingale~${(\FT{\Dir{\boldalpha_k}})}_k$. For every continuous compactly supported real-valued~$f$, it satisfies
\begin{align*}
\FT{\DF_\sigma}(f)\eqdef \int_{\msP(X)} \diff\DF_\sigma(\eta)\, e^{\imu \scalar{\eta}{f}}= \sum_{n=0}^\infty\frac{\imu^n}{n!}\,  Z_n(\sigma f^1, \dotsc,\sigma f^n)\comm
\end{align*}
where~$\imu=\sqrt{-1}$ is the imaginary unit,~$Z_n$ is the cycle index polynomial~\eqref{eq:MultiBell} of the~$n^\text{th}$ symmetric group and~$f^j$ denotes the~$j^{\textrm{th}}$ power of~$f$.

Furthermore, the map~$\sigma\mapsto \DF_\sigma$ is continuous with respect to the narrow topologies.
\end{thm}

The characteristic functional representation is new. It provides ---~in the unified framework of Fourier analysis~---
($a$) a new (although non-explicit) construction of~$\DF_{\sigma}$ as the unique probability measure on~$\msP(X)$ satisfying~$\FT{\DF_{\sigma}}=\klim \FT{\Dir{\boldalpha_k}}$ (see Cor.~\ref{c:BocMin}. Following~\cite{Ver07}, we call this construction a \emph{weak Fourier limit});
($b$) new proofs of known results on the tightness and asymptotics of families of Dirichlet--Ferguson measures (see Cor.s~\ref{c:Seth} and~\ref{p:Asympt}), proved, elsewhere in the literature, with ad hoc techniques;
($c$) the continuity statement in the Theorem, which strengthens~\cite[Thm.~3.2]{SetTiw81} concerned with norm-to-narrow continuity. This last result is sharp, in the sense that the domain topology cannot be relaxed to the vague topology.

\paragraph{Representations of~$SL_2$-currents and Bayesian non-parametrics} The Dirichlet--Ferguson measure~$\DF$, the \emph{gamma measure}~$\GP$~\cite{KonDaSStrUs98,TsiVerYor01} and the \emph{`multiplicative infinite-dimensional Lebesgue measure'}~$\mcL^+$~\cite{TsiVerYor01,Ver07} play an important r\^ole in a longstanding program~\cite{TsiVer03,TsiVerYor01,KonLytVer15} for the study of representations of measurable $SL_2$-current groups, i.e. spaces of $SL_2$-valued bounded measurable functions on a smooth manifold~$X$.
Within such framework, connections between these measures and Lie structures of \emph{special linear type} are not entirely surprising.
In particular, the measure~$\mcL^+$ is constructed (see~\cite[\S4.1]{Ver07}) as the weak Fourier limit for~$k\rar\infty$ of rescaled Haar measures on 
the identity connected components~$dSL^+_{k+1}$ in maximal toral ---~\emph{commutative}~--- subgroups of the special linear groups~$SL_{k+1}(\R)$. (For details on this construction see~\S\ref{ss:Algebra} below.) 

\smallskip

Relying on connections between cycle index polynomials and \emph{P\'olya Enumeration Theory}, we identify the special linear object acting on the Dirichlet distribution~$\Dir{\boldalpha_k}$ as the \emph{dynamical symmetry algebra}, in the sense of~\cite{Mil72,Mil73b}, of the Fourier transform~$\FT{\Dir{\boldalpha_k}}$.
In contrast with the case of~$\mcL^+$, we are able to detail the action of the whole ---~\emph{non-commutative}~--- dynamical symmetry algebra, and provide a suitable interpretation of this action in terms of Bayesian statistics. Indeed, one remarkable property~\cite{Fer73,Set94, NgTiaTan11} of Dirichlet measures is that their posterior distributions given knowledge on the occurrences of some categorical random variables are themselves Dirichlet measures with different parameters; that is, Dirichlet measures are \emph{self-conjugate priors}. We show how this property is related to the action of the dynamical symmetry algebra.
More precisely, for~$\boldalpha\in \Delta^{k-1}$ and~$\mbfp\in (\Z_0^+)^k$ denote by~$\Dir{\boldalpha}^\mbfp$ the posterior distribution of the prior~$\Dir{\boldalpha}$ given atoms of mass~$p_i$ at point~$i\in [k]$ (see property~\ref{i:3} in~\S\ref{ss:Dir}). We prove the following.

\begin{thm}[see Thm.~\ref{c:Due}]\label{t:DueIntro}
The dynamical symmetry algebra~$\mfg_k$ of the function~$\FT{\Dir{\boldalpha}}$ (see Def.~\ref{d:DSA}) is (isomorphic to) the Lie algebra~$\mfsl_{k+1}(\R)$ of real square matrices with vanishing trace. 
Furthermore, if~$\boldalpha$ is chosen in the interior of~$\Delta^{k-1}$, the universal enveloping algebra~$\mfU(\mfg_k)$ naturally acts on an infinite-dimensional linear space~$\mcO_{\Lambda_\boldalpha}$ detailed in the proof. 
Special subalgebras of~$\mfU(\mfg_k)$ may be identified, whose actions fix the linear span~$\mcO_{H_\boldalpha}\subset \mcO_{\Lambda_\boldalpha}$ of the family of characteristic functionals~$\{\FT{\Dir{\boldalpha}^\mbfp}\}_\mbfp$ varying~$\mbfp \in (\Z_0^+)^k$, or the linear span~$\mcO_{\Lambda_\boldalpha^+}\subset \mcO_{\Lambda_\boldalpha}$ of characteristic functionals of some distinguished improper priors of Dirichlet-categorical posteriors.
\end{thm}

Theorem~\ref{t:Uno} allows for a partial extension of this result to the infinite-dimensional case of~$\DF_\sigma$. Since~$\mcL^+$ and~$\mcG$ may be expressed as product measures with~$\DF$ as the only truly infinite-dimensional factor, cf.~\cite{TsiVerYor00}, we expect Theorem~\ref{t:DueIntro} to provide further algebraic insights on these measures.

\paragraph{Quasi-invariance of~$\DF$} (Quasi-)invariance properties of~$\DF$,~$\GP$ and~$\mcL^+$ have been studied with respect to different group actions~\cite{vReYorZam08,vReStu09,KonLytVer15,TsiVerYor01}.
Given~$(X,\sigma)$ a Riemannian manifold with normalized volume measure~$\sigma$, let~$G$ be some subgroup of (bi-)measurable isomorphisms of~$(X,\mcB(X))$.
We are interested in the quasi-invariance of~$\DF_\sigma$ with respect to the group action~$\psi.\eta\eqdef \psi_\pfwd \eta$ where $\psi$ is in~$G$,~$\eta$ is in~$\msP(X)$ and~$\psi_\pfwd \eta\eqdef \eta\circ\psi^{-1}$ denotes the push-forward of~$\eta$ via~$\psi$. When~$X=\mbbS^1$ and $G=\Diffeo(X)$, the quasi-invariance of~$\DF_\sigma$ with respect to a similar action was a key tool in the construction of stochastic dynamics on~$\msP(X)$ with~$\DF_\sigma$ or the related \emph{entropic measure}~$\mbbP_\sigma$ as invariant measures, see~\cite{Sha11,vReStu09}.

Whereas Theorem~\ref{t:Uno} allows for Bochner--Minlos and L\'evy Continuity related results to come into play, the non-multiplicativity of~$\FT{\DF_\sigma}$ (corresponding to the non-infinite-divisibility of the measure) immediately rules out the usual approach to quasi-invariance via Fourier transforms~\cite{AlbKonRoe98,TsiVerYor00,KonLytVer15, TsiVerYor01}. Other approaches to this problem  rely on finite-dimensional approximation techniques, variously concerned with approximating the space~\cite{vReStu09, vReYorZam08}, the~$\sigma$-algebra~\cite{KonLytVer15} or the acting group~\cite{Gor17, Ver07}. The common denominator here is for the approximation to be a filtration (cf. e.g.~\cite[Def.~9]{KonLytVer15}) ---~in order to allow for some kind of martingale convergence~--- and, possibly, for the approximating objects to be (embedded in) linear structures~(cf. e.g.~\cite{Ver07,vReStu09}).

The goal of Theorem~\ref{t:DueIntro} is ultimately to provide approximating sequences ---~at the same time of the space~$X$, the $\sigma$-algebra on~$\msP(X)$ and the acting group~--- that are suitable in the sense above.

\paragraph{Plan of the work}
Preliminary results are collected in~\S\ref{s:Prelim}, together with the definition and properties of Dirichlet measures and an account of the discretization procedure that we dwell upon in the following. In~\S\ref{s:t:Uno} we prove Theorem~\ref{t:Uno}. As a consequence, by the classical theory of characteristic functionals we recover known asymptotic expressions for~$\DF_{\beta\sigma}$ when~$\beta\rar0,\infty$ is a real parameter (Cor.~\ref{p:Asympt}, cf.~\cite[p.~311]{SetTiw81}), propose a \emph{Gibbsean} interpretation thereof (Rem.~\ref{r:Gibbsean}), and prove analogous expressions for the \emph{entropic measure}~$\mbbP^\beta_\sigma$ on compact Riemannian manifolds~\cite{Stu11}, generalizing the case~$X=\mbbS^1$~\cite[Prop.~3.14]{vReStu09}.
In the process of deriving Theorem~\ref{t:Uno} we obtain a moment formula for the Dirichlet distribution in terms of the cycle index polynomials~$Z_n$ (Thm.~\ref{t:MomDir}). In light of \emph{P\'olya Enumeration Theory} we interpret this result by means of a coloring problem~(\S\ref{ss:Color}). This motivates the study of the dynamical symmetry algebra~$\mfl_k$ of the Humbert function~${}_k\Phi_2$ resulting in the proof of Theorem~\ref{t:DueIntro}. Finally, in~\S\ref{ss:InfDimTDue} we study the limiting action of the dynamical symmetry algebra~$\mfl_k$ when~$k$ tends to infinity.

Some preliminary results in topology and measure theory are collected in the Appendix.

\section{Definitions and preliminaries}\label{s:Prelim}

\paragraph{Notation} Denote by~$\imu$ the imaginary unit, by~$\GF[a_n](t)$ (resp. by~$\EGF[a_n](t)$) the (exponential) generating function of the sequence~$\seq{a_n}_n$ of complex numbers, computed in the variable~$t$.

Let~$i,k,n$ be positive integers and set for~$1\leq i\leq k$ (the position of an element in a vector is stressed by a \emph{left} subscript)
\begin{align*}
\mbfy\eqdef& \seq{y_1,\dotsc, y_k} & \mbfe_i\eqdef&\seq{{}_10,\dotsc,0,{}_i1,0,\dotsc,{}_k0}\\
\uno\eqdef& \seq{{}_11,\dotsc,{}_k1} & \mbfy_{\hat \imath}\eqdef& \seq{y_1,\dotsc, y_{i-1}, y_{i+1}, \dotsc, y_k} \\
\vec{\mbfk}\eqdef& \seq{1,\dotsc,k} & \length{\mbfy}\eqdef& y_1+\cdots+y_k\fstop
%
%
\intertext{\indent
Write~$\mbfy>\zero$ for~$y_1,\dotsc,y_k>0$ and analogously for~$\mbfy\geq \zero$. 
Further set
}
[k]\eqdef& \set{1,\dotsc, k} & \pi\in\mfS_k\eqdef& \set{\textrm{bijections of }[k]}\\
\mbfy_\pi\eqdef& \seq{y_{\pi(1)}, \dotsc, y_{\pi(k)}} & \mbfy\compo\mbfz\eqdef& (y_1z_1,\dotsc, y_kz_k)\\
\mbfy^{\compo n}\eqdef& \underbrace{\mbfy\compo\dots\compo \mbfy}_{n \textrm{ times}} & \mbfy\cdot\mbfz\eqdef & y_1z_1+\cdots + y_kz_k \comm
\end{align*}
where~$\compo$ denotes the~\emph{Hadamard product} and we write~$\mbfy^{\compo\mbfz}=\tseq{y_1^{z_1},\dotsc, y_k^{z_k}}$ vs.~$\mbfy^\mbfz= y_1^{z_1}\cdots y_k^{z_k}$.
Given any $k$-variate complex-valued function~$f$, write
\begin{align*}
f(\mbfy)\eqdef f(y_1)\cdots f(y_k)&&& f^\compo(\mbfy)\eqdef \seq{f(y_1),\dotsc, f(y_k)} \fstop
\end{align*}

Finally, denote by $\Gamma$ the \emph{Euler Gamma function}, by $\Poch{\alpha}{k}\eqdef \Gamma(\alpha+k)/\Gamma(k)$ the \emph{Pochhammer symbol} of~\mbox{$\alpha\not\in \Z^-_0$}, by
$\Beta(y,z)\eqdef {\Gamma(y)\Gamma(z)}/{\Gamma(y+z)}$, resp.~$\Beta(\mbfy)\eqdef {\Gamma(\mbfy)}/{\Gamma(\length{\mbfy})}$,
the \emph{Euler Beta function}, resp. its multivariate analogue. 

\subsection{Combinatorial preliminaries}

\paragraph{Set and integer partitions} For a subset~$L\subset [n]$ denote by~$\tilde L$ the ordered tuple of elements in~$L$ in the usual order of~$[n]$. An \emph{ordered set partition} of $[n]$ is an ordered tuple~$\tilde\mbfL\eqdef(\tilde L_1,\tilde L_2\dotsc)$ of tuples~$\tilde L_i$ such that the corresponding sets~$L_i$, termed \emph{clusters} or blocks, satisfy~$\emp\subsetneq L_i\subset [n]$ and~$\sqcup_i L_i=[n]$. The order of the tuples in~$\tilde\mbfL$ is assumed ascending with respect to the cardinalities of the corresponding subsets and, subordinately, ascending with respect to the first element in each tuple. A \emph{set partition}~$\mbfL$ of~$[n]$ is the family of subsets corresponding to an ordered set partition. This correspondence is bijective.
For any set partition write~$\mbfL\vdash [n]$ and $\mbfL \vdash_r [n]$ if $\#\mbfL=r$, i.e. if $\mbfL$ has $r$ clusters.
A (\emph{integer}) \emph{partition~$\boldlambda$ of~$n$ into~$r$ parts} (write:~$\boldlambda\vdash_r n$) is an integer solution~$\boldlambda\geq \zero$ of the system,~$\vec\mbfn\cdot\boldlambda=n$,~$\length{\boldlambda}=r$; if the second equality is dropped we term $\boldlambda$ a (\emph{integer}) \emph{partition of $n$} (write:~$\boldlambda\vdash n$). We always regard a partition in its \emph{frequency representation}, i.e. as the tuple of its ordered frequencies (cf. e.g.~\cite[\S1.1]{And76}).
To a set partition \mbox{$\mbfL\vdash_r [n]$} one can associate in a unique way a partition~$\boldlambda(\mbfL)\vdash_r n$ by setting 
$\lambda_i(\mbfL)\eqdef \#\set{h\mid \#L_h=i}$.

\paragraph{Permutations and cycle index} A permutation $\pi$ in~$\mfS_n$ 
 is said to have \emph{cycle structure~$\boldlambda$}, write~$\boldlambda=\boldlambda(\pi)$, if~$\lambda_i$ equals the number of cycles in~$\pi$ of length~$i$ for each~$i$. Let~\mbox{$\mfS_{n}(\boldlambda)\subset \mfS_n$} be the set of permutations with~cycle structure~$\boldlambda$, so that~$\mfS_{n}(\boldlambda(\pi))=K_\pi$ the conjugacy class of~$\pi$ and~$\#\mfS_{n}(\boldlambda)=\Mnom{\boldlambda}\eqdef {n!}/({\boldlambda!\, \vec{\mbfn}^{\boldlambda}})$~\mbox{\cite[Prop.~I.1.3.2]{Sta01}}.

Let now $G<\mfS^n$ be any permutation group. The \emph{cycle index polynomial} of~$G$ is defined by
\begin{align*}
Z^G(\mbft)\eqdef \frac{1}{\# G} \sum_{\pi\in G} \mbft^{\boldlambda(\pi)}	\comm \mbft\eqdef \seq{t_1,\dotsc,t_n}\fstop
\end{align*}

We write $Z_n\eqdef Z^{\mfS_n}$ for the cycle index polynomial of~$\mfS_n$, satisfying, for~$\mbft\eqdef \seq{t_1,\dotsc,t_n}$ and $\mbft_k\eqdef\seq{t_1,\dotsc,t_k}$ with $k\leq n$, the identities
\begin{align}\label{eq:MultiBell}
Z_n(\mbft)= \frac{1}{n!}\sum_{\boldlambda\vdash n} \Mnom{\boldlambda} \mbft^\boldlambda \comm \qquad Z_n((a\uno)^{\compo\vec\mbfn}\compo\mbft)=a^n Z_n(\mbft) \quad a\in \R\fstop
\end{align}
and the recurrence relation
\begin{align}\label{eq:RecBellZ}
Z_{n}(\mbft)=\frac{1}{n}\sum_{k=0}^{n-1} Z_{k}(\mbft_k) \, x_{n-k}\comm Z_0(\emp)\eqdef 1\fstop 
\end{align}

\subsection{The Dirichlet distribution}\label{ss:Dir}
Denote the \emph{standard}, resp. \emph{corner}, $(k-1)$-dimensional simplex by
\begin{align*}
\Delta^{k-1}\eqdef\tset{\mbfy\in\R^k\mid \mbfy\geq \zero,\, \length{\mbfy}=1} \comm \Delta^{k-1}_*\eqdef \tset{\mbfz\in \R^{k-1}\mid \mbfz\geq \zero,\, \length{\mbfz}\leq 1} \fstop
\end{align*}

\begin{defs}[Dirichlet distribution] We denote by~$\Dir{\boldalpha}(\mbfy)$ the \emph{Dirichlet distribution with parameter~$\boldalpha\in \R_+^k$} (e.g.~\cite{NgTiaTan11}), i.e. the probability measure with density
\begin{align}\label{eq:DensityD}
\car_{\Delta^{k-1}}(\mbfy)\, \frac{\mbfy^{\boldalpha-\uno}}{\Beta(\boldalpha)}
\end{align}
with respect to the~$k$-dimensional Lebesgue measure on the hyperplane of equation~$\length{\mbfy}=1$ in~$\R^{k}$, concentrated on (the interior of)~$\Delta^{k-1}$. Alternatively, for any measurable~$A\subset \R^{k-1}$,
\begin{align*}
\Dir{\boldalpha}(A)=\int_{\Delta^{k-1}_*} \car_A(\mbfz) \prod_{i=1}^k z_i^{\alpha_i-1} \diff \mbfz \quad \textrm{where} \quad z_k\eqdef 1-\length{\mbfz} \fstop
\end{align*}
\end{defs}

Whereas both descriptions are common in the literature, the first one makes more apparent property~\ref{i:2} below.
Namely, write~`$\sim$' for `\emph{distributed as}' and let~$\mbfY$ be any $\Delta^{k-1}$-valued random vector.
The following properties of the Dirichlet distribution are well-known:
\begin{enumerate}[i.]
\item \emph{aggregation} (e.g.~\cite[p.~211, property $\textrm{i}^\circ$]{Fer73}). For $i=2,\dotsc, k$ set~$\mbfy_{\contraplus i}\eqdef(\mbfy+y_{i}\,\mbfe_{i-1})_{\hat{\imath}}$. Then,
\begin{align}\label{eq:Aggregation}
\mbfY\sim \Dir{\boldalpha}\implies& \mbfY_{\contraplus i}\sim \Dir{\boldalpha_{\contraplus i}} \fstop
\end{align}

\item\label{i:2} \emph{quasi-exchangeability} (or \emph{symmetry}). For all~$\pi\in \mfS_{k}$
\begin{align}\label{eq:Symmetry}
\mbfY\sim \Dir{\boldalpha}\implies \mbfY_\pi\sim\Dir{\boldalpha_\pi} \fstop
\end{align}

\item\label{i:3} \emph{Bayesian property} (e.g.~\cite[p.~212, property $\textrm{iii}^\circ$]{Fer73} for the case~$r=1$). Let~$\mbfW\in [k]^r$ be a vector of $[k]$-valued random variables and~$\mbfP\in (\Z_0^+)^k$ be the vector of occurrences defined by~$P_i\eqdef\#\set{j\in[r]\mid W_j=i}$. For~$\mbfp\in (\Z_0^+)^k$ let~$\mbfY$ be such that~$\Prob\set{P_i=p_i\mid \mbfY}=Y_i$ for all~$i\in [k]$ and denote by~$\Dir{\boldalpha}^\mbfp$ the distribution of~$\mbfY$ given~$\mbfP=\mbfp$, termed here the \emph{posterior distribution of~$\Dir{\boldalpha}$ given atoms with masses~$p_i$ at points~$i\in[k]$}. Then,
\begin{align}
\mbfY\sim \Dir{\boldalpha}\implies \Dir{\boldalpha}^\mbfp = \Dir{\boldalpha+\mbfp} \fstop
\end{align}
\end{enumerate}

\smallskip

Most properties of the Dirichlet distribution may be inferred from its characteristic functional~${}_k\Phi_2$, a confluent form of the \emph{$k$-variate Lauricella hypergeometric function}~${}_kF_D$ 
(see e.g.~\cite{Ext76}).
Recall the following representations of ${}_kF_D$~\cite[\S2.1]{Ext76} with~$\mbfb,\mbfs\in \C^k$, $a\in\C$ and~$c\in\C\setminus\Z_0^-$
\begin{align*}
{}_kF_D[a,\mbfb;c;\mbfs]\eqdef& \sum_{\mbfm\in \N_0^k} \frac{\Poch{a}{\length{\mbfm}} \Poch{\mbfb}{\mbfm} \mbfs^\mbfm}{ \Poch{c}{\length{\mbfm}} \mbfm! } & \norm{\mbfs}_\infty<1\\
 %
=& \frac{1}{\Beta(a,c-a)} \int_0^1 t^{a-1}(1-t)^{c-a-1}(\uno-t\mbfs)^{-\mbfb} \diff t & \Re c>\Re a>0
\end{align*}
and its \emph{confluent} form (or \emph{second $k$-variate Humbert function}~\cite[\emph{ibid.}]{Ext76}),~$\mbfb,\mbfs\in\C^k$
\begin{align}\label{eq:ConflLaur}
{}_k\Phi_2[\mbfb;c;\mbfs]\eqdef \lim_{\eps\rar 0^+} {}_kF_D[1/\eps;\mbfb;c;\eps\mbfs]=\sum_{\mbfm\in \N_0^k} \frac{\Poch{\mbfb}{\mbfm} \mbfs^\mbfm}{ \Poch{c}{\length{\mbfm}} \mbfm! } && c\in\C\setminus \Z_0^- \fstop
\end{align}

The distribution~$\Dir{\boldalpha}$ is moment determinate for any~$\boldalpha>\zero$ by compactness of $\Delta^{k-1}$.
Its moments are straightforwardly computed via the multinomial theorem as
\begin{equation}\label{eq:MomDir}
\begin{aligned}
\mu'_n[\mbfs,\boldalpha]\eqdef \int_{\Delta^{k-1}} (\mbfs\cdot \mbfy)^n \diff\Dir{\boldalpha}(\mbfy) 
=\sum_{\substack{\mbfm\in\N^k_0\\ \length{\mbfm}=n}} \tbinom{n}{\mbfm} \mbfs^\mbfm \tfrac{\Beta(\boldalpha+\mbfm)}{\Beta(\boldalpha)}=\frac{n!}{\Poch{\length{\boldalpha}}{n}} \sum_{\substack{\mbfm\in\N^k_0\\ \length{\mbfm}=n}} \frac{\mbfs^{\mbfm}}{\mbfm!} \Poch{\boldalpha}{\mbfm} \comm
\end{aligned}
\end{equation}
so that the characteristic functional of the distribution indeed satisfies (cf.~\cite[\S7.4.3]{Ext76})
\begin{equation}\label{eq:LapDir}
\begin{aligned}
\FT{\Dir{\boldalpha}}(\mbfs)\eqdef&\int_{\Delta^{k-1}} \exp(\imu \mbfs\cdot\mbfy)\diff\Dir{\boldalpha}(\mbfy) 
= 
\sum_{\mbfm\in\N^k_0}\frac{\Poch{\boldalpha}{\mbfm}}{ \mbfm! }\frac{\imu^{\length{\mbfm}} \mbfs^\mbfm}{\Poch{\length{\boldalpha}}{\length{\mbfm}}} \defeq {}_k\Phi_2[\boldalpha;\length{\boldalpha};\imu \,\mbfs] \fstop
\end{aligned}
\end{equation}

\subsection{The Dirichlet--Ferguson measure}
\paragraph{Notation} Everywhere in the following let~$(X,\tau(X))$ be a \emph{second countable locally compact Hausdorff} topological space with Borel $\sigma$-algebra~$\mcB$. We denote respectively by~$\cl A$,~$\interior A$,~$\bd A$ the \emph{closure}, \emph{interior} and \emph{boundary} of a set~$A\subset X$ with respect to~$\tau$. Recall (Prop.~\ref{p:Alex}) that any space~$(X,\tau(X))$ as above is Polish, i.e. there exists a metric~$d$, metrising~$\tau$, such that~$(X,d)$ is separable and complete; we denote by~$\diam A$ the \emph{diameter} of~$A\subset X$ with respect to any such metric~$d$ (apparent from context and thus omitted in the notation). 

Denote by~$\mcC_c(X)$ (resp.~$\mcC_b(X)$) the space of continuous compactly supported (resp. continuous bounded) functions on~$(X,\tau(X))$, (both) endowed with the topology of uniform convergence; by~$\mcC_0(X)$ the completion of~$\mcC_c(X)$, i.e. the space of continuous functions on~$X$ vanishing at infinity; by~$\Mb(X)$ (resp.~$\Mbp(X)$) the space of \emph{finite}, signed (resp. non-negative) \emph{Radon} measures on~$(X,\mcB(X))$ ---~the topological dual of~$\mcC_c(X)$ and~$\mcC_0(X)$~--- endowed with the the vague topology~$\tau_v(\Mb(X))$, i.e. the weak* topology, and the induced Borel~$\sigma$-algebra. Denote further by~$\msP(X)\subset \Mbp(X)$ (cf. Cor.~\ref{c:Vak}) the space of probability measures on~$(X,\mcB(X))$.
If not otherwise stated, we assume~$\msP(X)$ to be endowed with the vague topology~$\tau_v(\msP(X))$ and $\sigma$-algebra~$\mcB_v(\msP(X))$.
On~$\Mbp(X)$ (resp. on~$\msP(X)$) we additionally consider the narrow topology~$\tau_n(\Mbp(X))$ (resp.~$\tau_n(\msP(X))$), i.e. the topology induced by duality with~$\mcC_b(X)$.

Finally, given any measure~$\nu\in\Mb(X)$ and any bounded measurable function~$g$ on~$(X,\mcB(X))$, denote by~$\nu g$ the expectation of~$g$ with respect to~$\nu$ and by~$g^*\colon\nu\mapsto \nu g$ the linear functional induced by~$g$ on~$\Mb(X)$ via integration.

The following statement is well-known. A proof is sketched to establish further notation.
\begin{prop}\label{p:Alex}
A topological space~$(X,\tau(X))$ is second countable locally compact Hausdorff if and only if it is locally compact Polish, i.e. such that~$\tau(X)$ is a locally compact separable completely metrizable topology on~$X$.
Moreover, if~$(X,\mcB(X))$ additionally admits a fully supported diffuse measure~$\nu$, then~$(X,\tau(X))$ is perfect, i.e. it has no isolated points.
\begin{proof}[Sketch of proof] Let~$(\upalpha X,\tau(\upalpha X))$ denote the Alexandrov compactification of~$(X,\tau(X))$ and $\upalpha\colon X\rar \upalpha X$ denote the associated embedding. Notice that~$\upalpha X$ is Hausdorff, for~$X$ is locally compact Hausdorff; hence~$\upalpha X$ is metrizable, for it is second countable compact Hausdorff, and separable, for it is second countable metrizable, thus Polish by compactness. Finally, recall that~$X$ is (homeomorphic via $\upalpha$ to) a $G_\delta$-set in~$\upalpha X$ and every~$G_\delta$-set in a Polish space is itself Polish. The converse and the statement on perfectness are trivial.
\end{proof} 
\end{prop}

\paragraph{Partitions} Fix~$\sigma\in\msP(X)$. We denote by~$\mfP_k(X)$ the family of measurable non-trivial $k$-partitions of~$(X,\mcB,\sigma)$, i.e. the set of tuples~$\mbfX\eqdef \seq{X_1,\dotsc, X_k}$ such that
\begin{align*}
X_i\in \mcB\comm \sigma X_i>0\comm X_i\cap X_j=\emp \quad i,j\in [k], i\neq j\comm \cup_{i\in[k]} X_i=X \fstop
\end{align*}

Given~$\mbfX\in\mfP_k(X)$ we say that it \emph{refines}~$A$ in~$\mcB$ if~$X_i\subset A$ whenever $X_i\cap A\neq \emp$, respectively that it is a \emph{continuity partition for~$\sigma$} if~$\sigma(\bd X_i)=0$ for all~$i\in[k]$. We denote by~$\mfP_k(A\subset X)$, resp.~$\mfP_k(X,\tau(X),\sigma)$ the family of all such partitions.
Given~$\mbfX_1\in\mfP_{k_1}(X)$ and~$\mbfX_2\in\mfP_{k_2}(X)$ with~$k_1<k_2$ we say that~$\mbfX_2$ \emph{refines}~$\mbfX_1$, write~$\mbfX_1\preceq \mbfX_2$, if for every~$i\in [k_2]$ there  exists~$j_i\in [k_1]$ such that~$X_{2,i}\subset X_{1,j_i}$.
A sequence~$\seq{\mbfX_h}_h$ of partitions~$\mbfX_h\in\mfP_{k_h}(X)$ is termed a \emph{monotone null-array} if~$\mbfX_{h+1}\preceq\mbfX_h$ and $\hlim \max_{i\in [k_h]} \diam X_{h,i}=0$ (recall that~$\diam X_{h,i}$ vanishes independently of the chosen metric on~$(X,\tau(X))$, cf.~\cite[\S2.1]{Kal83}). We denote the family of all such null-arrays by~$\mfN\mfa(X)$. Analogously to partitions, we write with obvious meaning of the notation~$\mfN\mfa(A\subset X)$ and~$\mfN\mfa(X,\tau(X),\sigma)$. If~$\sigma$ is \emph{diffuse} (i.e. atomless), then~$\hlim \sigma X_{h,i_h}=0$ for every choice of~$X_{h,i_h}\in\mbfX_h$ with~$\seq{\mbfX_h}_h\in\mfN\mfa(X)$.

Given a (real-valued) simple function~$f$ and a partition~$\mbfX\in\mfP_k(X)$, we say that~$f$ is \emph{locally constant on~$\mbfX$ with values~$\mbfs$} if~$f\restr_{X_i}\equiv s_i$ constantly for every~$X_i\in \mbfX$.
Given a function~$f$ in~$\mcC_c$ we say that a sequence of (measurable) simple functions~$\seq{f_h}_h$ is a \emph{good approximation} of~$f$ if~$\abs{f_h}\uparrow_h \abs{f}$ and~$\hlim f_h=f$ pointwise. The existence of good approximations is standard.

\paragraph{The Dirichlet--Ferguson measure}
By a \emph{random probability} over~$(X,\mcB(X))$ we mean any probability measure on~$\msP(X)$.
For~$\mbfX\in\mfP_k(X)$ and~$\eta$ in~$\msP(X)$ set~$\eta^\compo\mbfX\eqdef \seq{\eta X_1,\dotsc, \eta X_k}$ and
\functionnn{\ev^\mbfX}{\msP(X)}{\Delta^{k-1}\subset \R^{k}}{\eta}{\eta^\compo\mbfX \fstop}

Recall (cf.~\cite{Sie22}) that, if~$\sigma\in\msP(X)$ is diffuse, then for every~$k\in\N_1$ and~$\mbfy\in \interior\Delta^{k-1}$ there exists~$\mbfX\in\mfP_k(X)$ such that~$\sigma^\compo\mbfX=\mbfy$.

\begin{defs}[Dirichlet--Ferguson measure] Fix~$\beta>0$ and~$\sigma\in\msP(X)$. The \emph{Dirichlet--Ferguson measure~$\DF_{\beta\sigma}$ with intensity~$\beta\sigma$}~\cite[\S1,~Def.~1]{Fer73} (or: \emph{Dirichlet}~\cite{LijPru10}, \emph{Poisson--Dirichlet}~\cite{Ver07}, \emph{Fleming--Viot with parent-independent mutation}~\cite{Fen10}; see e.g.~\cite[\S2]{Set94} for an explicit construction) is the unique random probability over~$(X,\mcB(X))$ such that
\begin{align}\label{eq:PartX}
\ev^{\mbfX}_\pfwd \DF_{\beta\sigma}=& \Dir{\beta \ev^{\mbfX}\sigma} \comm \mbfX\in\mfP_k(X)\comm k\in \N_1
\end{align}
(recall that~$\sigma^\compo\mbfX>\zero$). More explicitly, for every bounded measurable function $u\colon \Delta^{k-1}\rar\R$
\begin{equation}\label{eq:d:DirFer}
\int_{\msP(X)} u(\eta^\compo\mbfX) \diff\DF_{\beta\sigma}(\eta)=\int_{\Delta^{k-1}} u(\mbfy) \diff \Dir{\beta \sigma^\compo\mbfX}(\mbfy) \fstop
\end{equation}
\end{defs}

Existence was originally proved in~\cite{Fer73} by means of Kolmogorov Extension Theorem (cf.~Fig.~\ref{fig:KC} below). A construction on spaces more general than in our assumptions is given in~\cite{Kin06}. Other characterizations are available (see e.g.~\cite{Set94}). Since~$X$ is Polish (Prop.~\ref{p:Alex}), in~\eqref{eq:d:DirFer} it is in fact sufficient to consider~$u$ continuous with $\abs{u}<1$ and, by the Portmanteau Theorem,~$\mbfX\in\mfP_k(X,\tau(X),\sigma)$ (cf.~e.g.~\cite[p.~15]{Stu11}).

Let~$P$ be a~$\msP(X)$-valued random field on a probability space~$(\Omega,\msF,\mbbP)$ and recall the following properties of~$\DF_\sigma$, to be compared with those of~$\Dir{\boldalpha}$,
\begin{enumerate}[i.]
\item\label{i:DF1} \emph{realization properties}: If~$P\sim\DF_{\beta\sigma}$, then $P(\omega)=\sum_{i\in I}\eta_i(\omega) \delta_{x_i(\omega)}$ is $\mbbP$-a.s. purely atomic~\cite[\S4,~Thm.~2]{Fer73}, with~$\supp P(\omega)=\supp\sigma$~\cite[\S3,~Prop.~1]{Fer73} or~\cite{Maj92}. In particular, if~$\sigma$ is diffuse and fully supported, then~$I$ is countable and~$\set{x_i}_i$ is $\mbbP$-a.e.~dense in~$X$.
The sequence~$\seq{\eta_i}_i$ is distributed~\cite{Gri88} according to the \emph{stick-breaking process}. In particular,~$\mbbE\eta_i=\beta^{i-1}/(1+\beta)^i$. The r.v.'s~$x_i$'s are i.i.d. (independent also of the~$\eta_i$'s~\cite{EthKur94}) and $\sigma$-distributed.

\item\label{i:DF2} \emph{$\sigma$-symmetry}: for every measurable $\sigma$-preserving map~$\psi\colon X\rar X$, i.e. such that~$\psi_\pfwd \sigma=\sigma$,
\begin{align}\label{eq:Invar}
P\sim \DF_\sigma \implies \psi_\pfwd P \sim \DF_\sigma
\end{align}
(consequence of~\cite[Lem.~9.0]{Kal83} together with~\eqref{eq:PartX} and the quasi-exchangeability of~$\Dir{\boldalpha}$). In particular,~$P^\compo \mbfX$ is distributed as a function of~$\sigma^\compo\mbfX$ for every~$\mbfX\in\mfP_k(X)$ for every~$k$.

\item\label{i:DF3} \emph{Bayesian property}~\cite[\S3, Thm.~1]{Fer73}: Let~$\mbfW\eqdef\seq{W_1,\dotsc, W_r}$ be a sample of size~$r$ from~$P$, conditionally i.i.d., and denote by~$\DF^\mbfW_\sigma$ the distribution of~$P$ given~$\mbfW$, termed the \emph{posterior distribution of~$\DF_\sigma$ given atoms~$\mbfW$}. Then,
\begin{align*}
P\sim \DF_\sigma \implies (P\mid \mbfW) \sim \DF_{\sigma+\sum_j^r \delta_{W_j}} \fstop
\end{align*}
\end{enumerate}

\paragraph{Discretizations} 
In order to consider finite-dimensional marginalizations of~$\DF_{\beta\sigma}$, we introduce the following discretization procedure (cf.~\cite{RegSaz00} for a similar construction).
Any partition~$\mbfX\in\mfP_k(X)$ induces a discretization of~$X$ to~$[k]$ by collapsing~$X_i\in \mbfX$ to an arbitrary point in~$X_i$, uniquely identified by its index~$i\in[k]$, i.e. via the map~$\pr^{\mbfX}\colon X\supset X_i\ni x\mapsto i\in[k]$. The finite $\sigma$-algebra~$\sigma_0(\mbfX)$ generated by~$\mbfX$ induces then a discretization of~$\msP(X)$ to the space~$\msP([k])$ via the mapping~$\mu\mapsto \sum_i \mu X_i \,\delta_i$. Since the latter space is in turn homeomorphic to the standard simplex~$\Delta^{k-1}$ via the mapping~$\sum_i y_i\delta_i\mapsto \mbfy$, every choice of~$\mbfX\in\mfP_k(X)$ induces a discretization of~$\msP(X)$ to~$\Delta^{k-1}$ via the resulting composition~$\ev^\mbfX=\pr^\mbfX_\pfwd$.
It is then precisely the content of~\eqref{eq:PartX} that any partition~$\mbfX$ as above induces a discretization of the tuple~$((X,\sigma),(\msP(X),\DF_{\beta\sigma}))$ to the tuple~$(([k],\boldalpha),(\Delta^{k-1},\Dir{\boldalpha}))$, where~$\boldalpha\eqdef \beta \ev^\mbfX\sigma$ is identified with the measure~$\sum_i \alpha_i\delta_i$ on~$[k]$ (cf. Fig.~\ref{fig:KC} below).

Going further in this fashion, the subgroup~$\mfS_\mbfX$ of bi-measurable isomorphisms~$\psi$ of~$(X,\mcB(X))$ respecting~$\mbfX$, i.e. such that~$\psi^\compo(\mbfX)\eqdef\seq{\psi(X_1),\dotsc, \psi(X_k)}=\mbfX$ up to reordering, is naturally isomorphic to the symmetric group~$\mfS_k$, the bi-measurable isomorphism group~$\Iso([k])$ of~$[k]$. The canonical action of~$\mfS_\mbfX$ on~$X$, corresponding to the canonical action of~$\mfS_k$ on~$[k]$, lifts to the action of~$\mfS_k$ on~$\Delta^{k-1}$ by permutation of its vertices, that is, to the action on~$\msP([k])$ defined by~$\pi.\mbfy\eqdef \pi_\pfwd \mbfy$ under the identification of~$\mbfy$ with the measure~$\sum_i y_i\delta_i$.

\section{Proof of Theorem~\ref{t:Uno} and accessory results}\label{s:t:Uno}

\subsection{Finite-dimensional statements}\label{s:fDThm1}
Thinking of~$\boldalpha$ as a measure on~$[k]$ as in~\S\ref{s:Prelim}, the aggregation property~\eqref{eq:Aggregation} may be given a measure-theoretical interpretation too.
Indeed with the same notation of~\S\ref{ss:Dir}, for~$i\in [k-1]$ let additionally~$\mfs^{i}\colon [k]\rar [k-1]$ denote the~$i^{\textrm{th}}$ codegeneracy map of~$[k]$, i.e. the unique weakly order preserving surjection such that~$\#(\mfs^i)^{-1}(i)=2$. Then, up to the usual identification of~$\Delta^{k-1}$ with~$\msP([k])$, it holds that~$\mfs^i_\pfwd \mbfy=\mbfy_{\contraplus i}$ and one has $\mfs^i_\pfwd\mbfY\sim \mbfY_{\contraplus i}$. Thus, choosing $\mbfY\sim \Dir{\boldalpha}$, the aggregation property reads
$(\mfs^i_\pfwd)_\pfwd \Dir{\boldalpha}=\Dir{\mfs^i_\pfwd \boldalpha}$.

The following result is a rather obvious generalization of the latter fact, obtained by substituting degeneracy maps with arbitrary maps. We provide a proof for completeness.

\begin{prop}[Mapping Theorem for~{$\Dir{\boldalpha}$}]\label{p:MapDir}
Fix~$\boldalpha\in\R_+^k$. Then, for every~$g\colon [k]\rar [k]$
\begin{align*}
(g_\pfwd)_\pfwd \Dir{\boldalpha}=\Dir{g_\pfwd \boldalpha} \fstop
\end{align*}

\begin{proof} 
Define the \emph{additive contraction}~$\mbfy_{\contraplus\boldlambda}$ of a vector~$\mbfy$ with respect to~$\boldlambda\vdash k$ as
\begin{equation}\label{eq:AddContract}
\begin{aligned}
\mbfy_{\contraplus\boldlambda} \eqdef& (\underbrace{y_1,\dotsc, y_{\lambda_1}}_{\lambda_1}, \underbrace{y_{\lambda_1+1}+y_{\lambda_1+2}, \dotsc, y_{\lambda_1+2\lambda_2-1}+y_{\lambda_1+2\lambda_2}}_{2\lambda_2},\dotsc, \\
&\quad \underbrace{y_{\vec\mbfk\cdot\boldlambda-k\lambda_k+1}+\dots+y_{\vec\mbfk\cdot\boldlambda-(k-1)\lambda_k}, \dotsc, y_{\vec\mbfk \cdot\boldlambda-\lambda_k+1}+\dots+y_{\vec\mbfk\cdot\boldlambda}}_{k\lambda_k}) \comm
\end{aligned}
\end{equation}
whence inductively applying \eqref{eq:Aggregation} to any $\Delta^{k-1}$-valued random variable~$\mbfY$ yields $\mbfY\sim \Dir{\boldalpha} \implies \mbfY_{\contraplus\boldlambda}\sim \Dir{\boldalpha_{\contraplus\boldlambda}}$ for $\boldlambda \vdash k$.
Combining the latter with the quasi-exchangeability~\eqref{eq:Symmetry}, 
$\Dir{\boldalpha}$ satisfies
\begin{align}\label{eq:AggrGen}
\mbfY\sim \Dir{\boldalpha} \implies (\mbfY_{\pi})_{\contraplus\boldlambda}\sim \Dir{(\boldalpha_\pi)_{\contraplus\boldlambda}} \qquad \pi\in \mfS_k\comm \boldlambda \vdash k \fstop
\end{align}

For $\boldlambda\vdash k$ set~$\lambda_0\eqdef 0$ and define the map~$\contra\boldlambda\colon[k]\rar[\abs{\boldlambda}]$ by
\begin{align*}
\contra\boldlambda\colon i\mapsto \lambda_{j-1}+\ceiling{i/j} \quad\text{if}\quad i\in \set{(j-1)\lambda_{j-1}+1,\dotsc, j\lambda_j}
\end{align*}
varying~$j$ in~$[k]$, where~$\ceiling{\alpha}$ denotes the ceiling of~$\alpha$. It is readily checked that~$(\contra\boldlambda\circ \pi)_\pfwd \boldalpha=(\boldalpha_\pi)_{\contraplus\boldlambda}$ for any~$\pi$ in~$\mfS_k$. The proof is completed by exhibiting, for fixed~$g\colon [k]\rar [k]$, the unique partition~\mbox{$\boldlambda_g\vdash k$} and some permutation~$\pi_g\in\mfS_k$ such that~$g=\contra\boldlambda_g\circ \pi_g$.
To this end set~$L_{g,(i)}\eqdef g^{-1}[i]$ and
\begin{itemize}
\item[] $\mbfL_g\eqdef \seq{L_{g,(1)},\dotsc, L_{g,(k)}}$, where it is understood that~$L_{g,(i)}$ is omitted if empty; 
\item[] $\tilde\mbfL_g\eqdef(\tilde L_{1,1}, \tilde L_{1,2}, \dotsc, \tilde L_{2,1}, \dotsc)$ the ordered set partition associated to~$\mbfL_g$, where 
\item[] $\tilde L_{j,r}\eqdef \seq{\ell_{j,r,1}, \dotsc, \ell_{j,r,j}}$ denotes the~$r^\text{th}$ tuple of cardinality~$j$ in~$\tilde \mbfL_g$;
\end{itemize}

moreover, varying~$j$ in~$[k]$ and~$r$ in~$\floor{k/\lambda_j}$, where~$\floor{\alpha}$ denotes the floor of~$\alpha$, define~$\pi$ in~$\mfS_k$~by
\begin{align*}
\pi\colon i\mapsto \ell_{j,\,r,\, (i-\lambda_{j-1}-1 \mod j) +1} \quad \text{if}\quad \begin{cases} i\in \set{(j-1)\lambda_{j-1}+1,\dotsc, j\lambda_j}\\
\ceiling{(i-\lambda_{j-1}-1)/\lambda_j}=r
\end{cases} \fstop
\end{align*}

Finally set~$\pi_g\eqdef \pi^{-1}$ and~$\boldlambda_g\eqdef \boldlambda(\mbfL_g)$.
\end{proof}
\end{prop}

\begin{rem}
Assuming the point of view of conditional expectations rather than that of marginalizations,~\eqref{eq:PartX} may be restated as 
\begin{align*}
\mbbE_{\DF_{\beta\sigma}}[\emparg\vert \sigma_0(\mbfX)]=\mbbE_{\Dir{\beta \sigma^\compo \mbfX}}[\emparg]\comm
\end{align*}
where~$\sigma_0(\mbfX)$ denotes as before the~$\sigma$-algebra generated by some partition~$\mbfX\in\mfP_k(X)$. The aggregation property~\eqref{eq:Aggregation} is but an instance of the tower property of conditional expectations, whereas its generalization~\eqref{eq:AggrGen} is a consequence of the $\sigma$-symmetry of~$\DF_\sigma$.
\end{rem}

\begin{thm}[Moments of~$\Dir{\boldalpha}$]\label{t:MomDir}
Fix~$\boldalpha>\zero$ and~$\mbfs\in \R^k$. Then, the following identity holds
\begin{align}\label{eq:l:MomDir}
\mu_n'[\mbfs,\boldalpha]=\frac{n!}{\Poch{\length{\boldalpha}}{n}} \sum_{\substack{\mbfm\in\N^k_0\\ \length{\mbfm}=n}} \frac{\mbfs^{\mbfm}}{\mbfm!} \Poch{\boldalpha}{\mbfm}=\frac{n!}{\Poch{\length{\boldalpha}}{n}} Z_n(\mbfs^{\compo 1}\cdot\boldalpha, \dotsc,\mbfs^{\compo n}\cdot \boldalpha) \defeq \zeta_n[\mbfs,\boldalpha] \fstop
\end{align}

\begin{proof}
Let
\begin{align*}
\tilde\mu_n[\mbfs,\boldalpha]\eqdef& \Poch{\length{\boldalpha}}{n} (n!)^{-1} \mu_n'[\mbfs,\boldalpha]\comm & \tilde\zeta_n[\mbfs,\boldalpha]\eqdef& \Poch{\length{\boldalpha}}{n} (n!)^{-1}\zeta_n[\mbfs,\boldalpha] \fstop
\end{align*}

The statement is equivalent to $\tilde\mu_n=\tilde\zeta_n$, which we prove in two steps.

\medskip

\emph{Step 1}. The following identity holds
\begin{align}\label{eq:Dirichlet0.0}
\tilde\mu_{n-1}[\mbfs,\boldalpha+\mbfe_\ell]=\sum_{h=1}^n s_\ell^{h-1} \tilde\mu_{n-h}[\mbfs,\boldalpha] \fstop
\end{align}

By induction on $n$ with trivial (i.e. $1=1$) base step $n=1$.
\emph{Inductive step.} Assume for every~$\boldalpha>\zero$ and~$\mbfs$ in~$\R^k$
\begin{align}\label{eq:Dirichlet0}
\tilde\mu_{n-2}[\mbfs,\boldalpha+\mbfe_\ell]=\sum_{h=1}^{n-1} s_\ell^{h-1} \tilde\mu_{n-1-h}[\mbfs,\boldalpha] \fstop
\end{align}

Let~$\partial_j\eqdef \partial_{s_j}$ and notice that
\begin{equation}\label{eq:Dirichlet1}
\begin{aligned}
\partial_j \tilde\mu_n[\mbfs,\boldalpha] =&\sum_{\substack{\mbfm\in \N_0^k\\ \length{\mbfm}=n}} \frac{m_j \, \mbfs^{\mbfm-\mbfe_j}}{\mbfm !} \Poch{\boldalpha}{\mbfm}=
\sum_{\substack{\mbfm\in \N_0^k\\ \length{\mbfm}=n}} \frac{\mbfs^{\mbfm-\mbfe_j}}{(\mbfm-\mbfe_j) !} \alpha_j\Poch{\boldalpha+\mbfe_j}{\mbfm-\mbfe_j}\\
=& \alpha_j \sum_{\substack{\mbfm\in \N_0^k\\ \length{\mbfm}=n-1}} \frac{\mbfs^{\mbfm}}{\mbfm !}\Poch{\boldalpha+\mbfe_j}{\mbfm}= \alpha_j \tilde\mu_{n-1}[\mbfs,\boldalpha+\mbfe_j] \fstop
\end{aligned}
\end{equation}

If $k\geq 2$, we can choose $j\neq \ell$. Applying \eqref{eq:Dirichlet1} to both sides of \eqref{eq:Dirichlet0.0} yields
\begin{align*}
\partial_j \tilde\mu_{n-1}[\mbfs,\boldalpha+\mbfe_\ell]=& \alpha_j \tilde\mu_{n-2}[\mbfs,\boldalpha+\mbfe_j+\mbfe_\ell]\\
\partial_j \sum_{h=1}^n s_\ell^{h-1} \tilde\mu_{n-h}[\mbfs,\boldalpha]=& \sum_{h=1}^n s_\ell^{h-1} \alpha_j \tilde\mu_{n-h-1}[\mbfs,\boldalpha+\mbfe_j]\\
=& \alpha_j\sum_{h=1}^{n-1} s_\ell^{h-1} \tilde\mu_{n-h-1}[\mbfs,\boldalpha+\mbfe_j] \comm
\end{align*}
where the latter equality holds by letting~$\tilde\mu_{-1}\eqdef 0$. Letting now $\boldalpha'\eqdef \boldalpha+\mbfe_j$ and applying the inductive hypothesis~\eqref{eq:Dirichlet0} with $\boldalpha'$ in place of~$\boldalpha$ yields
\begin{align*}
\partial_j \tonde{\tilde\mu_{n-1}[\mbfs,\boldalpha+\mbfe_\ell]- \sum_{h=1}^n s_\ell^{h-1}\tilde\mu_{n-h}[\mbfs,\boldalpha]}=0
\end{align*}
for every $j\neq \ell$. By arbitrariness of $j\neq \ell$, the bracketed quantity is a polynomial in the sole variables~$s_\ell$ and~$\boldalpha$ of degree at most $n-1$ (obviously, the same holds also in the case $k=1$). As a consequence (or trivially if $k=1$), every monomial not in the sole variable~$s_\ell$ cancels out by arbitrariness of~$\mbfs$, yielding
\begin{align*}
\tilde\mu_{n-1}[\mbfs,\boldalpha+\mbfe_\ell]- \sum_{h=1}^n s_\ell^{h-1}\tilde\mu_{n-h}[\mbfs,\boldalpha]= \frac{s_\ell^{n-1}\Poch{\alpha_\ell+1}{n-1}}{(n-1)!}-\sum_{h=1}^n s_\ell^{h-1} \frac{s_\ell^{n-h}}{(n-h)!} \Poch{\alpha_\ell}{n-h} \fstop
\end{align*}

The latter quantity is proved to vanish as soon as
\begin{align*}
\frac{\Poch{\alpha+1}{n-1}}{(n-1)!}=\sum_{h=1}^n \frac{\Poch{\alpha}{n-h}}{(n-h)!} \comm  \textrm{or equivalently} \qquad \Poch{\alpha+1}{n-1}=\sum_{h=0}^{n-1} \frac{\Poch{\alpha}{h} (n-1)!}{h!}\comm
\end{align*}
in fact a particular case of the well-known \emph{Chu--Vandermonde identity} 
\begin{align*}
\Poch{\alpha+\beta}{n}=\sum_{k=0}^n \tbinom{n}{k}\Poch{\alpha}{k}\Poch{\beta}{n-k} \fstop
\end{align*}

\emph{Step 2.} It holds that $\tilde\mu_n=\tilde\zeta_n$. By strong induction on~$n$ with trivial (i.e. $1=1$) base step~$n=0$. \emph{Inductive step.}
Assume for every~$\boldalpha>\zero$ and~$\mbfs$ in~$\R^k$ that $\tilde\mu_{n-1}[\mbfs,\boldalpha]=\tilde\zeta_{n-1}[\mbfs,\boldalpha]$. Then
\begin{align*}
\partial_j \tilde\zeta_n[\mbfs,\boldalpha]=& \sum_{\boldlambda \vdash n} \frac{\Mnom{\boldlambda}}{n!} \sum_{h=1}^n \frac{\partial_j (\mbfs^{\compo h}\cdot \boldalpha)^{\lambda_h}}{(\mbfs^{\compo h}\cdot \boldalpha)^{\lambda_h} } \prod_{i=1}^n (\mbfs^{\compo i}\cdot \boldalpha)^{\lambda_i} \\
=&\sum_{\boldlambda \vdash n} \frac{\Mnom{\boldlambda}}{n!} \sum_{h=1}^n \frac{h\lambda_h s_j^{h-1}\alpha_j}{\mbfs^{\compo h}\cdot \boldalpha} \prod_{i=1}^n (\mbfs^{\compo i}\cdot \boldalpha)^{\lambda_i}\\
=&\alpha_j \sum_{h=1}^n s_j^{h-1} \sum_{\boldlambda \vdash n} \frac{h\lambda_h}{1^{\lambda_1}\lambda_1! \dotsc h^{\lambda_h}\lambda_h! \dotsc n^{\lambda_n}\lambda_n!} \frac{1}{\mbfs^{\compo h}\cdot \boldalpha} \prod_{i=1}^n (\mbfs^{\compo i}\cdot \boldalpha)^{\lambda_i}
\\
=&\alpha_j \sum_{h=1}^n s_j^{h-1} \sum_{\boldlambda \vdash n-h} \frac{\Mnom{\boldlambda}}{(n-h)!} \prod_{i=1}^{n-h} (\mbfs^{\compo i}\cdot\boldalpha)^{\lambda_i} \\
=&\alpha_j \sum_{h=1}^n s_j^{h-1} \tilde\zeta_{n-h}[\mbfs,\boldalpha] \fstop
\end{align*}

The inductive hypothesis, \eqref{eq:Dirichlet0.0} and \eqref{eq:Dirichlet1} yield
\begin{align*}
\partial_j \tilde\zeta_n[\mbfs,\boldalpha]=\alpha_j \sum_{h=1}^n s_j^{h-1} \tilde\mu_{n-h}[\mbfs, \boldalpha]=\partial_j \tilde\mu_n[\mbfs,\boldalpha] \fstop
\end{align*}

By arbitrariness of $j$ this implies that $\tilde\zeta_n[\mbfs,\boldalpha]-\tilde\mu_n[\mbfs,\boldalpha]$ is constant as a function of $\mbfs$ (for fixed~$\boldalpha$), hence vanishing by choosing~$\mbfs=\zero$.
\end{proof}
\end{thm}

\begin{rem}
Here, we gave an elementary combinatorial proof of the moment formula for~$\Dir{\boldalpha}$, independently of any property of the distribution. Notice for further purposes that, \emph{defining}~$\mu'_n[\mbfs,\boldalpha]$ as in~\eqref{eq:l:MomDir}, the statement holds with identical proof for all~$\boldalpha$ in~$\C^k$ such that~$\length{\boldalpha}\not\in\Z^-_0$.
For further representations of the moments see Remark~\ref{r:MomentsDF} below.
\end{rem}

\begin{prop}\label{p:LaurEGF}
The function ${}_k\Phi_2[t\mbfs;1;\boldalpha]$ is the exponential generating function of the polynomials~$Z_n$, in the sense that, for all~$\boldalpha\in \Delta^{k-1}$,
\begin{align*}
{}_k\Phi_2[\boldalpha; 1; t\mbfs]=& \EGF\quadre{Z_n\tonde{\mbfs^{\compo 1}\cdot\boldalpha, \dotsc, \mbfs^{\compo n}\cdot \boldalpha}}(t) \qquad \mbfs\in \R^k\comm t\in \R \fstop
\end{align*}

More generally,
\begin{align*}
{}_k\Phi_2[\boldalpha;\length{\boldalpha};t\mbfs]=\EGF\quadre{\frac{n!}{\Poch{\length{\boldalpha}}{n}}\, Z_n\tonde{\mbfs^{\compo 1}\cdot\boldalpha, \dotsc, \mbfs^{\compo n}\cdot \boldalpha}}(t)  \qquad \mbfs\in \R^k\comm t\in \R\fstop
\end{align*}

\begin{proof}
Recalling that ${}_k\Phi_2[\boldalpha;\length{\boldalpha};\mbfs]=\FT{\Dir{\boldalpha}}(\mbfs)$ by~\eqref{eq:LapDir} and noticing that $\length{\boldalpha}=1$, Theorem~\ref{t:MomDir} provides an exponential series representation for the characteristic functional of the Dirichlet distribution in terms of the cycle index polynomials of symmetric groups, viz.
\begin{align*}
\FT{\Dir{\boldalpha}}(\mbfs)=&\sum_{n=0}^\infty \frac{1}{n!} \, Z_n\tonde{(\imu\mbfs)^{\compo 1}\cdot\boldalpha, \dotsc,(\imu\mbfs)^{\compo n}\cdot\boldalpha}
\fstop
\end{align*}

Replacing~$\mbfs$ with~$-\imu t\mbfs$ above and using~\eqref{eq:MultiBell} to extract the term~$ t^n$ from each summand, the conclusion follows. The second statement has a similar proof.
\end{proof}
\end{prop}

\begin{rem}\label{r:TotPos}
It is well-known that the characteristic functional of a measure~$\mu$ on~$\R^d$ (or, more generally, on a nuclear space) is always \emph{positive definite}, i.e. it holds that
\begin{align*}
\forall n\in\N_0\quad \forall \mbfs_1,\dotsc, \mbfs_n\in \R^d \quad \forall \xi_1,\dotsc,\xi_n\in\C \qquad \sum_{h,k=1}^n \FT{\mu}(\mbfs_h-\mbfs_k) \, \xi_h \bar \xi_k\geq 0\comm
\end{align*}
where~$\bar\xi$ denotes the complex conjugate of~$\xi\in\C$.
Thus, whenever~$\boldalpha\in\R_+^k$, the functional~$\mbfs\mapsto {}_k\Phi_2[\boldalpha;\length{\boldalpha};\mbfs]$ is positive definite by~\eqref{eq:LapDir}. 
\end{rem}

The following Lemma also appeared in~\cite{LetPic18}.

\begin{lem}\label{l:Asympt}
There exist the narrow limits
\begin{align*}
\lim_{\beta\rar 0^+} \Dir{\beta\boldalpha}=& \length{\boldalpha}^{-1} \sum_{i=1}^{k} \alpha_i\delta_{\mbfe_i} & \text{and} &&
\lim_{\beta\rar +\infty} \Dir{\beta\boldalpha}=&\delta_{\length{\boldalpha}^{-1}\boldalpha} \fstop
\end{align*} 

\begin{proof}
Since $\Dir{\boldalpha}$ is moment determinate, it suffices ---~by compactness of $\Delta^{k-1}$ and  Stone--Weierstra{\ss} Theorem~--- to show the convergence of its moments. By Theorem~\ref{t:MomDir} (cf. also~\eqref{eq:MultiBell}), 
\begin{align*}
\mu_n'[\mbfs,\beta\boldalpha]\eqdef& \frac{n!}{\Poch{\beta\length{\boldalpha}}{n}} Z_n\tonde{\beta\mbfs^{\compo 1}\cdot\boldalpha, \dotsc, \beta \mbfs^{\compo n}\cdot\boldalpha}
=\frac{1}{\Poch{\beta\length{\boldalpha}}{n}} \sum_{r}^n \sum_{\boldlambda \vdash_r n} \Mnom{\boldlambda} \prod_i^n (\beta\, \mbfs^{\compo i}\cdot \boldalpha)^{\lambda_i}\\
=&\frac{1}{\Poch{\beta \length{\boldalpha}}{n}}\sum_r^n \sum_{\boldlambda \vdash_r n} \Mnom{\boldlambda} \beta^{\abs{\boldlambda}} \prod_i^n (\mbfs^{\compo i}\cdot \boldalpha)^{\lambda_i}
=\frac{1}{\Poch{\beta \length{\boldalpha}}{n}}\sum_r^n \beta^r \sum_{\boldlambda \vdash_r n} \Mnom{\boldlambda} \prod_i^n (\mbfs^{\compo i}\cdot \boldalpha)^{\lambda_i} \\
\underset{\beta\ll 1}{\approx}& \, \frac{1}{\beta \length{\boldalpha}\, \Gamma(n)}\, \beta \Mnom{\mbfe_n} (\mbfs^{\compo 1}\cdot \boldalpha)^{\abs{n\mbfe_1}}=\length{\boldalpha}^{-1}\boldalpha\cdot \mbfs^{\compo n} \comm\\
%
%
\underset{\beta\gg 1}{\approx} &\, \frac{1}{\beta^n \length{\boldalpha}^n}\, \beta^n \Mnom{n\mbfe_1} (\mbfs^{\compo 1}\cdot \boldalpha)^{\abs{n\mbfe_1}}= \length{\boldalpha}^{-n}(\mbfs\cdot \boldalpha)^n \fstop \qedhere
\end{align*}
\end{proof}
\end{lem}

As a consequence of the Lemma further confluent forms of~${}_k\Phi_2$ may be computed:

\begin{cor}[Confluent forms of~${}_k\Phi_2$]\label{c:Confluent}
There exist the limits
\begin{align*}
\lim_{\beta\rar 0^+} {}_k\Phi_2[\beta\boldalpha;\length{\beta\boldalpha};\mbfs]=&\length{\boldalpha}^{-1} \boldalpha\cdot \exp^\compo(\mbfs) \comm &
\lim_{\beta\rar +\infty} {}_k\Phi_2[\beta\boldalpha;\length{\beta\boldalpha};\mbfs]=& \exp(\length{\boldalpha}^{-1}\boldalpha\cdot\mbfs) \fstop
\end{align*}
\end{cor}

\subsection{Infinite-dimensional statements}
Together with the introductory discussion, Proposition~\ref{p:MapDir} suggests the following Mapping Theorem for~$\DF_\sigma$, to be compared with the analogous result for the Poisson random measure~$\PP_\sigma$ over $(X,\mcB(X))$ (see e.g.~\cite[\S2.3 and \emph{passim}]{Kin93}). The~$\sigma$-symmetry of~$\DF_{\beta\sigma}$ and the quasi-exchangeability and aggregation property of~$\Dir{\boldalpha}$ are trivially recovered from the Theorem by~\eqref{eq:PartX}.


\begin{thm}[Mapping theorem for $\DF_\sigma$]\label{t:Mapping}
Let~$(X,\tau(X),\mcB(X))$ and~$(X',\tau(X'),\mcB(X'))$ be second countable locally compact Hausdorff spaces,~$\nu$ a non-negative finite measure on~$(X,\mcB(X))$ and $f\colon (X,\mcB(X))\rar (X',\mcB'(X))$ be any measurable map. Then,
\begin{align*}
(f_\pfwd)_\pfwd\DF_\nu=\DF_{f_\pfwd\nu} \fstop
\end{align*}

\begin{proof}
Choosing $\mbfX\eqdef (g^{-1}(1), \dotsc, g^{-1}(k))$, the characterization~\eqref{eq:d:DirFer} is equivalent to the requirement that
$(g_\pfwd)_\pfwd \DF_\nu= \Dir{g_\pfwd \nu}$
for any~$g\colon X\rar [k]$ such that every~$\nu$-representative of~$g$ is surjective, which makes~$\mbfX$ non-trivial for~$\nu$. Denote by~$\mcS(X, \nu,k)$ the family of such functions and notice that if $h\in \mcS(X',f_\pfwd\nu,k)$, then~$g\eqdef h\circ f\in\mcS(X,\nu,k)$. The proof is now merely typographical:
\begin{align*}
(h_\pfwd)_\pfwd (f_\pfwd)_\pfwd \DF_\nu= (g_\pfwd)_\pfwd \DF_\nu = \Dir{g_\pfwd \nu}= \Dir{h_\pfwd (f_\pfwd \nu)} \comm
\end{align*}
where the second equality suffices to establish that $(f_\pfwd)_\pfwd \DF_\nu$ is a Dirichlet--Ferguson measure by arbitrariness of $h$, while the third one characterizes its intensity as $f_\pfwd \nu$.
\end{proof}
\end{thm}

We denote by~$\msP(\msP(X))$ the space of probability measures on~$(\msP(X),\mcB_n(\msP(X)))$, endowed with the narrow topology~$\tau_n(\msP(\msP(X)))$ induced by duality with~$\mcC_b(\msP(X))$. We are now able to prove the following more general version of Theorem~\ref{t:Uno}.

\begin{thm}[Characteristic functional of~$\DF_{\beta\sigma}$]\label{t:UnoBis}
Let~$(X,\tau(X),\mcB(X))$ be a second countable locally compact Hausdorff space,~$\sigma$ a probability measure on~$X$ and fix~$\beta>0$. Then, 
\begin{align}\label{eq:LapDF}
\forall f\in \mcC_c \qquad \FT{\DF_{\beta\sigma}}(tf^*)= \EGF\quadre{n!\Poch{\beta}{n}^{-1} Z_n\tonde{\beta \sigma f^1, \dotsc,\beta\sigma f^n}}(\imu\,t) \comm t\in\R \fstop
\end{align}

Moreover, the map~$\nu\mapsto \DF_{\nu}$ is narrowly continuous on~$\Mbp(X)$.

\begin{proof}
\emph{Characteristic functional.} Fix~$f$ in~$\mcC_c$ and let~$\seq{f_h}_{h}$ be a good approximation of~$f$, locally constant on~$\mbfX_h\eqdef\seq{X_{h,1},\dotsc,X_{h,k_h}}$ with values~$\mbfs_h$ for some~$\seq{\mbfX_h}_h\in\mfN\mfa(X)$. Fix~$n>0$ and set~$\boldalpha_h\eqdef \beta\sigma^\compo\mbfX_h$. Choosing~$u\colon \Delta^{k_h-1}\rar\R$, $u\colon\mbfy\mapsto (\mbfs_h\cdot \mbfy)^n$ in~\eqref{eq:d:DirFer} yields
\begin{align*}
\mu^{\prime\, \DF_{\beta\sigma}}_n[f_h^*]\eqdef& \int_{\msP(X)} (f_h^*\eta)^n\diff\DF_{\beta\sigma}(\eta)=\int_{\Delta^{k_h-1}} (\mbfs_h\cdot \mbfy)^n \diff \Dir{\beta\ev^{\mbfX_h}\sigma}(\mbfy)=\mu^{\prime}_n[\mbfs_h,\boldalpha_h]\comm
\intertext{hence, by Theorem~\ref{t:MomDir},}
\mu^{\prime\, \DF_{\beta\sigma}}_n[f_h^*]=&n!\Poch{\beta}{n}^{-1} Z_n\tonde{\mbfs_h^{\compo 1}\cdot\boldalpha_h, \dotsc, \mbfs_h^{\compo n}\cdot\boldalpha_h}=n!\Poch{\beta}{n}^{-1} Z_n\tonde{\beta\sigma f_h^1, \dotsc, \beta\sigma f_h^n} \comm
\end{align*}
thus, by Dominated Convergence Theorem, continuity of~$Z_n$ and arbitrariness of~$f$,
\begin{align*}
\forall f\in \mcC_c \qquad \mu^{\prime\, \DF_{\beta\sigma}}_n[tf^*]=&n!\Poch{\beta}{n}^{-1} Z_n\tonde{t^1 \beta\sigma f^1, \dotsc, t^n \beta\sigma f^n} \comm t\in \R \fstop
\end{align*}

Using~\eqref{eq:MultiBell} to extract the term~$t^n$ from~$Z_n$ and substituting~$t$ with~$\imu \, t$ on the right-hand side, the conclusion follows by definition of exponential generating function.

\medskip

\emph{Continuity.} Assume first that~$(X,\tau(X))$ is compact. By compactness of~$(X,\tau(X))$, the narrow and vague topology on~$\msP(X)$ coincide and~$\msP(X)$ is compact as well by Prokhorov Theorem. Let~$\seq{\nu_h}_{h\in \N}$ be a sequence of finite non-negative measures narrowly convergent to~$\nu_\infty$. Again by Prokhorov Theorem and by compactness of~$\msP(X)$ there exists some $\tau_n(\msP(\msP(X)))$-cluster point~$\DF_\infty$ for the family $\set{\DF_{\nu_h}}_{h}$.
By narrow convergence of~$\nu_h$ to~$\nu_\infty$, continuity of~$Z_n$ and absolute convergence of~$\FT{\DF_{\emparg}}(f)$, it follows that~$\hlim \FT{\DF_{\nu_h}}=\FT{\DF_{\nu_\infty}}$ pointwise on~$\mcC_c(X)$, hence, by Corollary~\ref{c:Vak}, it must be~$\DF_\infty=\DF_{\nu_\infty}$.

In the case when~$X$ is not compact, recall the notation established in Proposition~\ref{p:Alex}, denote by~$\mcB(\upalpha X)$ the Borel $\sigma$-algebra of~$(\upalpha X,\tau(\upalpha X))$ and by~$\msP(\upalpha X)$ the space of probability measures on~$(\upalpha X,\mcB(\upalpha X))$. 
By the Continuous Mapping Theorem there exists the narrow limit~$\tau_n(\msP(X))$-$\hlim \upalpha_\pfwd\nu_h=\upalpha_\pfwd\nu_\infty$, thus, by the result in the compact case applied to the space~$(\upalpha X,\mcB_\upalpha)$ together with the sequence~$\upalpha_\pfwd \nu_h$,
\begin{align}\label{eq:t:Continuity}
\tau_n(\msP(\msP(X)))\textrm{-}\!\hlim \DF_{\upalpha_\pfwd\nu_h}= &\DF_{\upalpha_\pfwd \nu_\infty} \fstop
\end{align}

The {narrow} convergence of~$\nu_h$ to~$\nu_\infty$ implies that~$\upalpha_\pfwd\nu_\infty$ does not charge the point at infinity in~$\upalpha X$, hence the measure spaces $(X,\mcB(X),\nu_*)$ and $(\upalpha X, \mcB(\upalpha X), \upalpha_\pfwd\nu_*)$ are isomorphic for~$*=h,\infty$ via the map~$\upalpha$, with inverse~$\upalpha^{-1}$ defined on~$\im\upalpha\subsetneq \upalpha X $. The continuity of~$\upalpha^{-1}$ and the Continuous Mapping Theorem together yield the narrow continuity of the map~$({\upalpha^{-1}}_\pfwd)_\pfwd$. The conclusion follows by applying~$({\upalpha^{-1}}_\pfwd)_\pfwd$ to~\eqref{eq:t:Continuity} and using the Mapping Theorem~\ref{t:Mapping}.
\end{proof}
\end{thm}

\begin{rem}\label{r:MomentsDF}
Different representations of the univariate moments of the Dirichlet--Ferguson measure have also appeared, without mention to~$Z_n$, in~\cite[Eq.~(17)]{Reg98} (in terms of incomplete Bell polynomials, solely in the case when~$X\Subset \R^+$ and $f=\id_{\R}$) and in~\cite[proof of Prop.~3.3]{LetPic18} (in implicit recursive form).
Representations of the multi-variate moments have also appeared in~\cite[Prop.~7.4]{KerTsi01} (in terms of summations over `color-respecting' permutations, in the case~$\beta=1$), in~\cite[(4.20)]{EthKur94} and~\cite[Lem.~5.2]{Fen10} (in terms of summations over constrained set partitions).
\end{rem}

In the case when~$\nu_h$ converges to~$\nu_\infty$ in total variation, the continuity statement in the Theorem and the asymptotics for~$\beta\rar0$ in Corollary~\ref{p:Asympt} below were first shown in~\cite[Thm.~3.2]{SetTiw81}, relying on Sethuraman's stick-breaking representation. The following result was also obtained, again with different methods, in~\cite{SetTiw81}.

\begin{cor}[Tightness of Dirichlet--Ferguson measures~{\cite[Thm.~3.1]{SetTiw81}}]\label{c:Seth} Under the same assumptions as in Theorem~\ref{t:UnoBis}, let~$M\subset \Mbp(X)$ be such that~$\overline M\eqdef \set{\nu/\nu X\mid \nu\in M}$ is a tight, resp. narrowly compact, family of finite non-negative measures. Then, the family~$\set{\DF_\nu}_{\nu\in M}$ is itself tight, resp. narrowly compact.
\end{cor}

\begin{cor}[Asymptotic expressions]\label{p:Asympt}
With the same assumptions as in Theorem~\ref{t:UnoBis} there exist for all~$f$ in~$\mcC_c$ and complex~$t$ the limits
\begin{align}
\label{eq:AsymptDFLap}
\lim_{\beta\downarrow 0} \FT{\DF_{\beta\sigma}}(tf^*)=&\sigma \exp(\imu\,t f) & \text{and} && \lim_{\beta\rar \infty} \FT{\DF_{\beta\sigma}}(tf^*)=& \exp(\imu\,t\, \sigma f)
\intertext{corresponding to the narrow limits}
\label{eq:AsymptDF}
\DF^0_\sigma\eqdef \lim_{\beta\downarrow 0} \DF_{\beta\sigma}=&\delta_\pfwd\sigma & \text{and} && \DF^\infty_\sigma\eqdef \lim_{\beta\rar \infty} \DF_{\beta\sigma}=&\delta_\sigma \comm
\end{align}
where, in the first case,~$\delta\colon X\rar \msP(X)$ denotes the Dirac embedding~$x\mapsto \delta_x$.

\begin{proof}
The existence of~$\DF^0_\sigma$ and~$\DF^\infty_\sigma$ as narrow cluster points for~$\set{\DF_{\beta\sigma}}_{\beta>0}$ follows by Corollary~\ref{c:Seth}. Retaining the notation established in Theorem~\ref{t:UnoBis}, Corollary~\ref{c:Confluent} yields for all~$k$
\begin{align}
\nonumber
\lim_{\beta\downarrow 0} \FT{\DF_{\beta\sigma}}(f_k^*)=& \sigma \exp(\imu f_k) & \text{and} && \lim_{\beta\rar\infty} \FT{\DF_{\beta\sigma}}(f_k^*)=& \exp(\imu\, \sigma f_k) \comm
\intertext{hence, by Dominated Converge,}
\label{eq:Thm1.1}
\klim \lim_{\beta\downarrow 0} \FT{\DF_{\beta\sigma}}(f_k^*)=& \sigma \exp(\imu f) & \text{and} && \klim\lim_{\beta\rar\infty} \FT{\DF_{\beta\sigma}}(f_k^*)=& \exp(\imu\, \sigma f) \fstop
\end{align}

Furthermore, recalling that~$\abs{f_k}\leq \abs{f}$ one has
\begin{align}\label{eq:Thm1.2}
\abs{\FT{\DF_{\beta\sigma}}(f^*)-\FT{\DF_{\beta\sigma}}(f_k^*)}\leq e^{\norm{f}} \int_{\msP(X)} \diff \DF_{\beta\sigma}(\eta) \abs{f-f_k}^*\eta \leq e^{\norm{f}}\norm{f-f_k}\comm
\end{align}
hence the order of the limits in each left-hand side of~\eqref{eq:Thm1.1} may be exchanged, for the convergence in~$k$ is uniform with respect to~$\beta$. This shows~\eqref{eq:AsymptDFLap}.
\end{proof}
\end{cor}

By Theorem~\ref{t:UnoBis},~$\beta\sigma$ may be substituted with any sequence~$\seq{\beta_h\sigma_h}_h$ with~$\hlim \beta_h=0,\infty$ and~$\set{\sigma_h}_h$ a tight family.
Observe that, despite the similarity with Lemma~\ref{l:Asympt}, Corollary~\ref{p:Asympt} is not a direct consequence of the former, since the evaluation map~$\ev^{\mbfX}$ is never continuous.

\begin{rem}[A Gibbsean interpretation]\label{r:Gibbsean}
Corollary~\ref{p:Asympt} states that, varying~$\beta\in[0,\infty]$, the map $\DF_{\beta\emparg}\colon \msP(X)\rar\msP(\msP(X))$ is a (continuous) interpolation between the two extremal maps~$\DF^0_\emparg=\delta^{(0)}_\pfwd$ and $\DF^\infty_\emparg=\delta^{(1)}$, where~$\delta^{(0)}\eqdef\delta\colon X\rar \msP(X)$ and~$\delta^{(1)}\eqdef \delta\colon \msP(X)\rar \msP(\msP(X))$.
These asymptotic distributions may be interpreted ---~at least formally~--- in the framework of statistical mechanics. In order to establish some lexicon, consider a physical system at \emph{inverse temperature}~$\beta$ driven by a \emph{Hamiltonian}~$H$.

Let~$Z_\beta^H\eqdef \av{\exp(-\beta H)}$,~$F_\beta\eqdef -\beta^{-1}\ln Z_\beta^H$ and~$G_\beta\eqdef (Z_\beta^H)^{-1}\exp(-\beta H)$ respectively denote the \emph{partition function}, the \emph{Helmholtz free energy} and (the distribution of) the \emph{Gibbs measure} of the system. It was heuristically argued in~\cite[\S3.1]{vReStu09} that ---~at least in the case when~$(X,\mcB,\sigma)$ is the unit interval~---
\begin{align*}
\diff \DF_{\beta\sigma}(\eta) = \frac{e^{-\beta \, S(\eta)}}{Z_\beta}\diff \DF^*_\sigma(\eta) \comm
\end{align*}
where:~$S$ is now an \emph{entropy functional} (rather than an energy functional),~$Z_\beta$ is a normalization constant and~$\beta$ plays the r\^ole of the {inverse temperature}. Here, $\DF^*_\sigma$ denotes a \emph{non-existing} (!) uniform distribution on~$\msP(X)$.
Borrowing again the terminology, this time in full generality, one can say that for small~$\beta$ (i.e. large \emph{temperature}), the system \emph{thermalizes} towards the ``uniform'' distribution $\delta_\pfwd\sigma$ induced by the reference measure $\sigma$ on the base space, while for large~$\beta$ it \emph{crystallizes} to $\delta_\sigma$, so that all randomness is lost.
Consistently with property~\ref{i:DF1} of~$\DF_\sigma$, we see that~$\mbbE_{\DF_\sigma^\infty}\eta_i=0$ and~$\mbbE_{\DF_\sigma^0}\eta_i=\delta_{i1}$ for all~$i$, where~$\delta_{ab}$ denotes the Kronecker symbol; in fact, both statements hold with probability 1.

It is worth noticing that a different interpretation for the parameter~$\beta$ has been given in~\cite{LetPic18}, where the latter is regarded as a `time' parameter in the definition of a PCOC.
\end{rem}

\begin{rem} By the Continuous Mapping Theorem, both the continuity statement in Theorem~\ref{t:UnoBis} and the asymptotic expressions in Corollary~\ref{p:Asympt} hold, mutatis mutandis, for every narrowly continuous image of~$\DF_{\beta\sigma}$, hence, for instance, for the \emph{entropic measure}~$\P^\beta_\sigma$~\cite{vReStu09,Stu11}. This generalizes~\cite[3.14]{vReStu09} and the discussion for the entropic measure thereafter.
\end{rem}

\begin{cor}[Alternative construction of~$\DF_{\beta\sigma}$]\label{c:BocMin}
Assume there exists a nuclear function space~$\mcS\subset \mcC_0(X)$, continuously embedded into~$\mcC_0(X)$ and such that~$\mcS\cap \mcC_c(X)$ is norm-dense in~$\mcC_0(X)$ and dense in~$\mcS$.
Then, there exists a unique Borel probability measure on the dual space~$\mcS'$, namely~$\DF_{\beta\sigma}$, whose characteristic functional is given by the extension of \eqref{eq:LapDF} to~$\mcS$.

\begin{proof} By the classical Bochner--Minlos Theorem (see e.g.~\cite[\S4.2, Thm.~2]{GelVil64}), it suffices to show that the extension to~$\mcS$, say~$\chi$, of the functional~\eqref{eq:LapDF} is a characteristic functional.
By the convention in~\eqref{eq:RecBellZ}, $\chi(\zero_\mcS)=\chi(\zero_{\mcC_c(X)})=1$. The (sequential) continuity of~$\chi$ on~$\mcS$ follows by that on~$\mcC_0(X)$ and the continuity of the embedding~$\mcS\subset \mcC_0(X)$. It remains to show the positivity (see Rmk.~\ref{r:TotPos}) of~$\chi$, which can be checked only on~$\mcS\cap \mcC_c(X)$ by $\norm{\emparg}$-density of the inclusions~$\mcS\cap\mcC_c(X)\subset \mcC_0(X)$. The positivity of~$\chi$ restricted to~$\mcC_c(X)$ follows from the positivity of~${}_k\Phi_2$ in Remark~\ref{r:TotPos} by approximation of~$f$ with simple functions as in the proof of Theorem~\ref{t:UnoBis}.
\end{proof}
\end{cor}

\begin{rem}
Let us notice that the assumption of Corollary~\ref{c:BocMin} is satisfied, whenever~$X$ is (additionally) either finite (trivially), or a differentiable manifold, or a topological group (by the main result in~\cite{Aar72}).
In particular, when~$X=\R^d$, we can choose~$\mcS=\mcS(\R^d)$, the space of Schwartz functions on~$\R^d$.
\end{rem}

\medskip

Consider the map~$\Giry\colon \msP(\msP(X))\rar \msP(X)$ defined by
\begin{align*}
(\Giry(\mu))A=\int_{\msP(X)} \diff \mu(\eta)\, \eta A \qquad A\in \mcB(X)\comm \mu\in\msP(\msP(X)) \fstop
\end{align*}

Since~$f^*$ is $\tau_n(\msP(X))$-continuous for every~$f\in\mcC_b(X)$ and bounded by~$\norm{f}$, the map~$\Giry$ is continuous.

\begin{cor} For fixed~$\beta\in(0,\infty)$, the map~$\DF_{\beta\emparg}\colon \msP(X)\rar\msP(\msP(X))$ is a homeomorphism onto its image, with inverse~$\Giry$.

\begin{proof} The continuity of~$\DF_{\beta\emparg}$ is proven in Theorem~\ref{t:UnoBis}. By e.g.~\cite[Thm.~3]{Fer73} for all~$f\in\mcC_c(X)$ one has~$\DF_{\beta\sigma}f^*=\sigma f$, hence~$\Giry$ inverts~$\DF_{\beta\emparg}$ on its image.
\end{proof}
\end{cor}

\begin{figure}[hbt!]
\small\center
\begin{tikzcd}[row sep=small, execute at end picture={
\foreach \Valor/\Nombre in   
{%
  tikz@f@1-1-1/a, tikz@f@1-2-2/b, tikz@f@1-5-2/c, tikz@f@1-4-1/d,%
  tikz@f@1-2-2/e, tikz@f@1-2-3/f, tikz@f@1-3-5/g, tikz@f@1-3-6/h,%
  tikz@f@1-6-5/i, tikz@f@1-6-6/j, tikz@f@1-5-2/k, tikz@f@1-5-3/l%
}
{
\coordinate (\Nombre) at (\Valor);
}
\fill[pattern=north east lines, pattern color=grey1,opacity=0.1]
  (a) -- (b) -- (c) -- (d) -- cycle;
\fill[pattern=north east lines, pattern color=grey1,opacity=0.5]
  (e) to (f) to [bend right=6] (g) to [bend left=12] (e);
\fill[pattern=north east lines, pattern color=grey1,opacity=0.5]
  (k) to (l) to [bend right=6] (i) to [bend left=11] (k);
\fill[pattern=north west lines, pattern color=grey1,opacity=0.3]
  (g) -- (h) -- (j) -- (i) -- cycle;
\fill[pattern=north west lines, pattern color=grey1,opacity=0.3]
  (e) -- (f) -- (l) -- (k) -- cycle;
\fill[pattern=north east lines, pattern color=grey1,opacity=0.1]
  (f) to [bend right=6] (g) -- (i) to [bend left=6] (l) -- cycle;
  }
]
\msP(\Delta^{k-1})
\\
\cdots \ar[no tail]{u} & \msP(\Delta^{k-1}) \ar{l} \ar["{g_\pfwd}_\pfwd" ']{ul} & \msP(\Delta^k) \ar[leftarrow, no head]{ddd}{\Dir{\beta\emparg}} \ar["{\mfs^i_\pfwd}_\pfwd"']{l} & \cdots \ar{l}
\\
&&&& \msP(\msP(X)) \ar[near start, bend left=5, "{\pr^{\mbfX_{k+1}}_\pfwd}_{\!\!\!\!\!\!\!\!\!\!\!\!\!\pfwd}" ']{ull} \ar[crossing over, near start, bend left=10, "{\pr^{\mbfX_k}_\pfwd}_{\!\!\!\!\!\!\pfwd}"]{ulll} \ar{r}{{f_\pfwd}_\pfwd} & \msP(\msP(X'))
\\
\Delta^{k-1} \ar[no head, near end]{uu}{\Dir{\beta\emparg}}
\\
\cdots \ar[no tail]{u} & \Delta^{k-1} \ar{l} \ar["g_\pfwd" ']{ul} \ar{uuu}{\Dir{\beta\emparg}} & \Delta^k \ar["\mfs^i_\pfwd" ']{l} \ar[leftarrow, no head]{dd}{\delta} & \cdots \ar{l}
\\
&&\phantom{a} && \msP(X) \ar{uuu}{\DF_{\beta\emparg}} \ar[crossing over, near start, bend left=10, "\pr^{\mbfX_k}_\pfwd"]{ulll} 	\ar[near start, bend left=5, "\pr^{\mbfX_{k+1}}_\pfwd" ']{ull} \ar{r}{f_\pfwd} & \msP(X') \ar{uuu}{\DF_{\beta\emparg}}
\\
{[k]} \ar[hook, no head, near end]{uu}{\delta} &&\phantom{a}&&&
\\
\cdots & {[k]} \ar{l} \ar["g" ']{ul} \ar[hook]{uuu}{\delta} & {[k+1]} \ar[no head, hook]{uu}{} \ar["\mfs^i" ']{l} 
& \cdots \ar{l}
\\
&&&& X \ar{r}{f} \ar[hook]{uuu}{\delta} \ar[near start, bend left=10, "\pr^{\mbfX_k}"]{ulll} 	\ar[near start, bend left=5, "\pr^{\mbfX_{k+1}}" ']{ull} & X' \ar[hook]{uuu}{\delta}
\end{tikzcd}
\captionsetup{singlelinecheck=off}
\caption[foo]{Many properties of Dirichlet(--Ferguson) measures can be phrased in terms of the commutation of some diagrams. The commutation of dashed squares of the diagram above, from left to right, respectively corresponds to
\begin{itemize}
\item the symmetry property~\eqref{eq:Symmetry} when~$g=\pi\in\mfS_k$ and, more generally, Proposition~\ref{p:MapDir};
\item the aggregation property~\eqref{eq:Aggregation};
\item the marginalization~\eqref{eq:PartX} (recall that~$\pr^{\mbfX}_\pfwd=\ev^\mbfX$);
\item the symmetry property~\eqref{eq:Invar} when~$f=\psi$ is measure preserving and, more generally, Theorem~\ref{t:Mapping};
\end{itemize}
the commutation of the solid sub-diagram delimited by the two dashed triangles corresponds to the requirement of Kolmogorov consistency.
}
\label{fig:KC}
\end{figure}

\section{Proof of Theorem~\ref{t:DueIntro} and accessory results}\label{s:t:Due}

\subsection{Finite-dimensional statements}\label{ss:Color}

\paragraph{Multisets} Given a set $S$, a (\emph{finite integer-valued}) \emph{$S$-multi-set} is any function $f\colon S\rar \N_1$ such that~$\#f$ is finite, where~$\#$ denotes integration on~$S$ with respect to the counting measure. We denote any such multiset by~$\bag{\mbfs_\boldalpha}$, where~$\mbfs\eqdef \seq{s_1,\dotsc, s_k}$ is $S$-valued with mutually different entries and~$\boldalpha\eqdef \seq{f(s_1),\dotsc, f(s_k)}$.
We term the set~$[\mbfs]\eqdef\set{s_1,\dotsc, s_k}$ the \emph{underlying set} to~$\bag{\mbfs_\boldalpha}$ and put 
\begin{align*}
[\mbfs_\boldalpha]\eqdef \set{(s_1,1), \dotsc, (s_1,\alpha_1), \dotsc, (s_k,1), \dotsc, (s_k,\alpha_k)} \fstop
\end{align*}

Recall that the number of $[n]$-multi-sets with cardinality~$r$ is~$\Poch{r}{n}/r!$~(see e.g.~\cite[\S{I.1.2}]{Sta01}). 

\subsubsection{A coloring problem}\label{sss:Color} An interpretation of the moments formula~\eqref{eq:l:MomDir} may be given in enumerative combinatorics, by means of \emph{P\'olya Enumeration Theory}~(PET, see e.g.~\cite{Pol37}); a minimal background is as follows.
Let~$G<\mfS_n$ be a permutation group acting on~$[n]$ and $[\mbfs]\eqdef\set{s_1,\dotsc,s_k}$ denote a set of (distinct) colors. A \emph{$k$-coloring} of~$[n]$ is any function~$f$ in~$[\mbfs]^{[n]}$, where we understand the elements~$s_1,\dotsc, s_k$ of~$[\mbfs]$ as placeholders for different {colors}. Whenever these are irrelevant, given a $k$-coloring~$f$ of~$[n]$ we denote by~$\tilde f$ the unique function in~$[k]^{[n]}$ such that~$s_{\tilde f(\emparg)}=f(\emparg)$.

We say that two $k$-colorings~$f_1,f_2$ of~$[n]$ are \emph{$G$-equivalent} if $f_1\circ \pi=f_2$ for all~$\pi$ in~$G$.

\begin{thm}[P\'olya~{\cite[\S4]{Pol37}}]\label{t:Polya}
Let~$G<\mfS_n$ be a permutation group acting on~$[n]$ and~$a_{h_1,\dotsc, h_k}$ be the number of $G$-inequivalent $k$-colorings of~$[n]$ into~$k$ colors with exactly~$h_i$ occurrences of the~$i^\text{th}$ color. Then, the (multivariate) generating function~$\GF[a_{h_1,\dotsc,h_k}](\mbft)$ satisfies
\begin{align}\label{eq:Polya}
\GF[a_{h_1,\dotsc,h_k}](\mbft)=Z^G\tonde{p_{k,1}[\mbft], \dotsc, p_{k,n}[\mbft]} \comm
\end{align}
where $p_{k,i}[\mbft]\eqdef \uno\cdot\, \mbft^{\compo i}$ with~$\uno\in \R^k$ denotes the~$i^\textrm{th}$ \emph{$k$-variate power sum symmetric polynomial}.
\end{thm}

In the following we consider an extension of PET to multisets of colors and explore its connections ---~arising in the case~$G=\mfS_n$~--- with the Dirichlet distribution~$\Dir{\boldalpha}$.
A different approach in terms of colorings, limited to the case~$\length{\boldalpha}=1$, was briefly sketched in~\cite[\S7]{KerTsi01}.

Let~$\bag{\mbfs_\boldalpha}$ be an integer-valued multiset with~$\boldalpha\in\R_+^k$, henceforth a~\emph{palette}. As before, we understand the elements~$s_1,\dotsc, s_k$ of its underlying set~$[\mbfs]$ as placeholders for different {colors}, and the elements~$(s_i,1), \dotsc, (s_i,\alpha_i)$ of~$[\mbfs_\boldalpha]$ as placeholders for different \emph{shades} of the same color~$s_i$. For a $k$-coloring~$f$ of~$[n]$ we say that~$\phi$ in~$[\mbfs_\boldalpha]^{[n]}$ is a \emph{shading} of~$f$ (and an~\emph{$\boldalpha$-shading} of~$[n]$) if~$\phi(\emparg)_1=f(\emparg)$, and that two shadings are $G$-equivalent if so are the corresponding colorings.

\begin{cor}[Counting shadings]\label{c:Shades}
Let~$G<\mfS_n$ be a permutation group acting on~$[n]$ and~$b_{h_1,\dotsc,h_k}^\boldalpha$ be the number of $G$-inequivalent $\boldalpha$-shadings of~$[n]$ with exactly~$h_i$ occurrences of the $i^\text{th}$~color. Then,
\begin{align*}
\GF[b_{h_1,\dotsc,h_k}^\boldalpha](\mbfs)=Z^G\tonde{\boldalpha\cdot\mbfs^{\compo 1}, \dotsc, \boldalpha\cdot\mbfs^{\compo n}} \comm \mbfs\in\R^k\fstop
\end{align*}

\begin{proof}
Set~$r_i\eqdef \length{\seq{\alpha_1,\dotsc,\alpha_i}}$ and~$r_0\eqdef 0$. For every~$r_k$-coloring~$g$ of~$[n]$ let
\begin{align*}
\phi_\boldalpha(x)\eqdef (s_i,\tilde g(x)-r_{i-1})\quad \text{if} \quad \tilde g(x)\in\set{r_{i-1}+1,\dotsc, r_i}
\end{align*}
varying~$i$ in~$[k]$. It is readily seen that, for every fixed~$\boldalpha$, this correspondence is bijective and preserves $G$-equivalence.
Thus the number~$a_{h_{1,1},\dotsc, h_{1,\alpha_1}, \dotsc, h_{k,1},\dotsc, h_{k,\alpha_k}}$ of $G$-inequivalent $r_k$-colorings of~$[n]$ with exactly~$h_{i,j}$ occurrences of the~$(r_{i-1}+j)^\text{th}$~color is also the number of~$G$-inequivalent \mbox{$\boldalpha$-shadings} of~$[n]$ with exactly~$h_{i,j}$ occurrences of the~$j^\text{th}$ shade of the $i^\text{th}$~color. By Theorem~\ref{t:Polya} this is the coefficient of the monomial
$t_1^{h_{1,1}}\cdots t_{r_1}^{h_{1,\alpha_1}} \cdots t_{r_{k-1}+1}^{h_{k,1}}\cdots t_{r_k}^{h_{k,\alpha_k}}$
in~$Z^G\tonde{\uno\cdot\, \mbft, \dotsc, \uno\cdot\, \mbft^{\compo n}}$ with~$\uno\in\R^{r_k}$.
By definition,
\begin{align*}
b_{h_1,\dotsc, b_{h_k}}^\boldalpha=\sum_{\substack{h_{1,1},\dotsc, h_{1,\alpha_1}, \dotsc, h_{k,1}, \dotsc, h_{k,\alpha_k}\\ \sum_j^{\alpha_i} h_{i,j}=h_i}} a_{h_{1,1},\dotsc, h_{1,\alpha_1}, \dotsc, h_{k,1},\dotsc, h_{k,\alpha_k}} \comm
\end{align*}
which equals the coefficient of the monomial~$s_1^{h_1}\dotsc s_k^{h_k}$ in
\begin{align*}
Z^G\tonde{\uno\cdot\, \mbft^{\compo 1}, \dotsc, \uno\cdot\, \mbft^{\compo n}}=&Z^G\tonde{\boldalpha\cdot \mbfs^{\compo 1},\dotsc, \boldalpha\cdot \mbfs^{\compo n}} \comm \mbft\eqdef(\underbrace{s_1,\dotsc, s_1}_{\alpha_1},\dotsc, \underbrace{s_k,\dotsc, s_k}_{\alpha_k}) \fstop \qedhere
\end{align*}
\end{proof}
\end{cor}

\begin{cor}\label{c:FracShades}
Let $\mcS_{n,k,r}$ denote the set of~$\mfS_n$-equivalence classes~$\phi^\bullet$ of $\boldalpha$-shadings of~$[n]$ such that~$\boldalpha\geq \zero$ and~$\length{\boldalpha}=r$. Then, the probability~$p_{h_1,\dotsc, h_k}^\boldalpha$ of some~$\phi^\bullet$ uniformly drawn from~$\mcS_{n,k,r}$ having exactly~$h_i$ occurrences of the $i^\text{th}$~color satisfies
\begin{align*}
\GF[p_{h_1,\dotsc,h_k}^\boldalpha](\mbfs)=\mu'_n[\mbfs,\boldalpha] \fstop
\end{align*}

\begin{proof} The number of palettes with total number of shades~$r$ equals the number~$\Poch{r}{n}/n!$ of integer-valued $[n]$-multisets of cardinality~$r$, thus, choosing~$r=\length{\boldalpha}$,
\begin{align*}
p_{h_1,\dotsc, h_k}^\boldalpha=n!\Poch{r}{n}^{-1} b_{h_1,\dotsc, h_k}^\boldalpha=n!\Poch{\length{\boldalpha}}{n}^{-1} b_{h_1,\dotsc, h_k}^\boldalpha \comm
\intertext{hence, by homogeneity}
\GF[p_{h_1,\dotsc, h_k}^\boldalpha](\mbfs)=n! \Poch{\length{\boldalpha}}{n}^{-1} \GF[b_{h_1,\dotsc, h_k}^\boldalpha](\mbfs) \fstop
\end{align*}

The conclusion follows by Corollary~\ref{c:Shades} and Theorem~\ref{t:MomDir}.
\end{proof}
\end{cor}

The study of~$\Dir{\boldalpha}$ in the case when~$\length{\boldalpha}=1$ is singled out as computationally easiest (as suggested by Theorem~\ref{t:MomDir}, noticing that~$\Poch{1}{n}=n!$), $\boldalpha$ representing in that case a \emph{probability} on~$[k]$, as detailed in~\S\ref{s:Prelim}. For these reasons, this is often the only case considered (cf. e.g.~\cite{KerTsi01}). On the other hand though, the general case when~$\boldalpha>\zero$ is the one relevant in Bayesian non-parametrics, since posterior distributions of Dirichlet-categorical and Dirichlet-multinomial priors do not have probability intensity. The above coloring problem suggests that the case when~$\boldalpha\in(\Z^+)^k$ is interesting from the point of view of PET, since it allows for some natural operations on palettes, corresponding to functionals of the distribution.

Indeed, we can change the number of colors and shades in a palette~$\bag{\mbfs_\boldalpha}$ by composing any permutation of the indices~$[k]$ with the following elementary operations: 
\begin{itemize}
\item $(i)$ `widen', respectively $(ii)$ `narrow the color spectrum', by adding a color, say~$s_{k+1}$, respectively removing a color, say~$s_k$. That is, we consider new palettes~$\bag{(\mbfs\oplus s_{k+1})_{\boldalpha\oplus\alpha_{k+1}}}$, respectively~$\bag{(s_1,\dotsc, s_{k-1})_{\tseq{\alpha_1,\dotsc,\alpha_{k-1}}}}$;
\item $(iii)$ `reduce color resolution' by regarding two different colors, say~$s_i$ and~$s_{i+1}$, as the same, relabeled~$s_i$. In so doing we regard the shades of the former colors as distinct shades of the new one, so that it has~$\alpha_i+\alpha_{i+1}$ shades. That is, we consider the new palette~$\bag{(\mbfs_{\hat \imath})_{\boldalpha_{\contraplus i}}}$;
\item $(iv)$ `enlarge', respectively $(v)$ `reduce the color depth', by adding a shade, say the~$\alpha_{i+1}^\text{th}$, to the color~$s_i$, respectively removing a shade, say the~$\alpha_i^\text{th}$, to the color~$s_i$. This latter operation we allow only if~$\alpha_i>1$, so to make it distinct from removing the color~$s_i$ from the palette. That is, we consider the new palettes~$\bag{\mbfs_{\boldalpha+\mbfe_i}}$, resp.~$\bag{\mbfs_{\boldalpha-\mbfe_i}}$ when $\alpha_i>1$.
\end{itemize}


Increasing the color resolution of a multi-shaded color, say~$s_k$ with $\alpha_k>1$~shades, by splitting it into two colors, say~$s_k'$ and~$s_{k+1}$ with~$\alpha_k'>0$ and $\alpha_{k+1}>0$ shades respectively and such that~$\alpha_k'+\alpha_{k+1}=\alpha_k$, is not an elementary operation. It can be obtained by widening the spectrum of the palette by adding a color~$s_{k+1}$ with~$\alpha_{k+1}$ shades and reducing the color depth of the color~$s_k$ to~$\alpha_k'$. Thus, this operation is not listed above. 
We do not allow for the number of shades of a color to be reduced to zero: although this is morally equivalent to removing that color, the latter operation amounts more rigorously to remove the color placeholder from the palette.

The said elementary operations are of two distinct kinds: $(i)$--$(iii)$ alter the number of colors in a palette, while $(iv)$--$(v)$ fix it. We restrict our attention to the latter ones and ask how the probability~$p_{h_1,\dotsc,h_k}^\boldalpha$ changes under them.
By Corollary~\ref{c:FracShades} this is equivalent to study the corresponding functionals of the~$n^\text{th}$ moment of the Dirichlet distribution. 
For fixed~$k$, we address all the moments at once, by studying the moment generating function
\begin{align*}
{}_k\Phi_2[\boldalpha;\length{\boldalpha}; t\mbfs]=\EGF[\GF[p_{h_1,\dotsc,h_k}^\boldalpha](\mbfs)](t) \fstop
\end{align*}

Namely, we look for natural transformations yielding the mappings
\begin{align}\label{eq:DSAOp}
E_{\pm i} \,\,{}_k\Phi_2 [\boldalpha;\length{\boldalpha}; \mbfs] = C_\boldalpha\,\, {}_k\Phi_2 [\boldalpha\pm \mbfe_i;\length{\boldalpha}\pm1; \mbfs] \comm
\end{align}
where~$C_\boldalpha$ is some constant, possibly dependent on~$\boldalpha$.
Here `natural' means that we only allow for meaningful \emph{linear} operations on generating functions: addition, scalar multiplication by variables or constants, differentiation and integration.
For practical reasons, it is convenient to consider the following construction.

\begin{defs}[Dynamical symmetry algebra of~${}_k\Phi_2$]\label{d:DSA}
Denote by~$\mfg_k$ the minimal Lie algebra containing the linear span of the operators~${E_{\pm 1}, \dotsc,  E_{\pm k}}$ in~\eqref{eq:DSAOp} endowed with the bracket induced by their composition.
Following~\cite{Mil72}, we term the Lie algebra~$\mfg_k$ the \emph{dynamical symmetry algebra} of the function~${}_k\Phi[\boldalpha;\mbfs]\eqdef {}_k\Phi_2[\boldalpha; \length{\boldalpha};\mbfs]$, characterized below.
\end{defs} 

\subsubsection{Dynamical symmetry algebras}\label{sss:DSA}
We compute now the dynamical symmetry algebra of the function~${}_k\Phi[\boldalpha;\mbfs]\eqdef {}_k\Phi_2[\boldalpha; \length{\boldalpha};\mbfs]$, in this section always regarded as the meromorphic extension~\eqref{eq:LapDir} of the Fourier transform of~$\FT{\Dir{\boldalpha}}(\mbfs)$ in the \emph{complex} variables~$\boldalpha,\mbfs\in \C^k$. The choice of complex variables is merely motivated by this identification and every result in the following concerned with complex Lie algebras holds verbatim for their split real form. 
For dynamical symmetry algebras of Lauricella hypergeometric functions see~\cite{Mil73,Mil72} and references therein; we refer to~\cite{Hum72} for the general theory of Lie algebra (representations) and for Weyl groups' theory.

\paragraph{Notation and definitions}
Denote by~$\mbfE_{i,j}$ varying $i,j\in[k+1]$ the canonical basis of~$\mathrm{Mat}_{k+1}(\C)$, with~$[\mbfE_{i,j}]_{m,n}=\delta_{mi}\,\delta_{nj}$, where~$\delta_{ab}$ is the Kronecker delta, and by~$A^*$ the conjugate transpose of a matrix~$A$. The following is standard.

\begin{lem}[A presentation of~$\mfsl_{k+1}(\C)$] 
For $i,j=0,\dotsc, k$ with~$j>i$ set
\begin{align*}
e_{i,j}\eqdef \mbfE_{i+1,j+1} \comm h_{i,j}\eqdef \mbfE_{i+1,i+1}-\mbfE_{j+1,j+1} \comm f_{j,i}\eqdef e_{i,j}^* \fstop
\end{align*}

Then, the complex Lie sub-algebra~$\mfl_k$ of~$\mfg\mfl_{k+1}(\C)$ generated by these vectors is~$\mfl_k=\mfsl_{k+1}(\C)$, with~$\mfsl_2$-triples
\begin{align*}
\set{e_{i,j},h_{i,j}, f_{j,i}}_{\substack{i=0,\dotsc, k\\j=i+1,\dotsc, k}} \fstop
\end{align*}

Denote further by~$\mff_k<\mfl_k$ the sub-algebra spanned by~$\set{e_{i,j}, f_{j,i}, h_{i,j}}_{i,j\in [k]}$. Then,~$\mff_k\cong\mfsl_k(\C)$.
\end{lem}

Everywhere in the following we regard~$\mfl_k$ together with the distinguished Cartan sub-algebra $\mfh_k<\mfl_k$~of diagonal traceless matrices spanned by the basis~$\set{h_{0,j}}_{j\in[k]}$; the root system~$\Psi_k$ induced by~$\mfh_k$, with simple roots~$\gamma_j$ corresponding to the~$\mfsl_2$-triples of the vectors~$e_{j-1,j}$ for~$j\in [k]$; positive, resp. negative, roots~$\Psi_k^\pm$ corresponding to the spaces of strictly upper, resp. strictly lower, triangular matrices~$\mfn^\pm_k$.
The inclusion~$\mff_k<\mfl_k$ induces the decomposition of vector spaces (\emph{not} of algebras)
\begin{align*}
\mfl_k=\mfr_k^-\oplus \mfh_1\oplus \mff_k\oplus \mfr_k^+\comm \textrm{where}\qquad \mfr_k^+\eqdef \C\set{e_{0,j}}_{j\in [k]}\comm \mfr_k^-\eqdef \C\set{f_{j,0}}_{j\in [k]}\comm \mfh_1= \C\set{h_{0,1}} \fstop
\end{align*}

%
The subscript~$k$ is omitted whenever apparent from context.

\medskip

For fixed~$\boldalpha\in \C^k$ regard~${}_k\Phi[\boldalpha;\emparg]$ as a formal power series and let~$f_{\boldalpha}\colon \C^{2k+1}_{\mbfs,\mbfu, t}\longrightarrow \C$ be
\begin{align}\label{eq:basis2}
f_\boldalpha= f_\boldalpha(\mbfs, \mbfu, t)\eqdef& {}_k\Phi[\boldalpha;\mbfs] \mbfu^\boldalpha t^{\length{\boldalpha}} \fstop
\end{align}

Let~$\Alpha\subset\C^k$. It is readily seen that the functions~$\set{f_\boldalpha}_{\boldalpha\in\Alpha}$ are (finitely) linearly independent, since so are the functions~$\set{f_\boldalpha(\uno,\mbfu,1)\propto\mbfu^\boldalpha}_{\boldalpha\in \Alpha}$. Set
\begin{align*}
\mcO_{\Alpha}\eqdef \bigoplus_{\boldalpha\in \Alpha} \C\!\set{f_\boldalpha} \comm \mcO_\boldalpha\eqdef \mcO_{\set{\boldalpha}} \comm \mcO\eqdef \mcO_{\C^k}
\end{align*}
and define the following differential operators, acting formally on~$\mcO$,
\begin{equation}\label{eq:DiffOpPhi}
\begin{aligned}
E_{\alpha_i}\eqdef& u_it(s_i\partial_{s_i}+u_i\partial_{u_i}-(\mbfs\cdot \nabla^\mbfs)\partial_{s_i})\comm & E_{\alpha_i,-\alpha_j}\eqdef& u_iu_j^{-1}\tonde{(u_i-u_j)\partial_{s_i}+u_i\partial_{u_i}} \,,\\
 E_{-\alpha_i}\eqdef& (u_i t)^{-1}(s_i-\mbfs\cdot\nabla^{\mbfs}-t\partial_t+1) \comm & J_{\alpha_i}\eqdef& t\partial_t+u_i\partial_{u_i}-1\,,\\
\end{aligned}
\end{equation}
where $i,j\in[k]$, $i\neq j$ and $\nabla^\mbfy\eqdef \seq{\partial_{y_1},\dotsc,\partial_{y_k}}$ for $\mbfy=\mbfu,\mbfs$. Term the operators~$E_{\alpha_i}$, resp.~$E_{-\alpha_i}$, \emph{raising}, resp. \emph{lowering}, \emph{operators}. Finally, let~$\mfg_k$ be the complex linear span of the operators~\eqref{eq:DiffOpPhi} endowed with the bracket induced by their composition.

\paragraph{Actions on spaces of holomorphic functions}
Let~$\Lambda_\boldalpha\eqdef \boldalpha+\Z^k$ and set, for every~$\ell\in\R^+$,
\begin{align*}
\Lambda_\boldalpha^+\eqdef \set{\boldalpha'\in \Lambda_\boldalpha\mid \length{\boldalpha'}>0}\comm H_\boldalpha\eqdef\set{\boldalpha'\in \Lambda_\boldalpha^+\mid \boldalpha'>\zero} \comm M_{\boldalpha,\ell}\eqdef \set{\boldalpha'\in \Lambda_\boldalpha^+\mid \length{\boldalpha'}=\ell} \fstop 
\end{align*}

Notice that if~$\Re^\compo\boldalpha>\zero$, the space~$\mcO_{\Lambda_\boldalpha^+}$ is a space of holomorphic functions~$\mcO(\C^k_\mbfs\times (\C\setminus \R^-_0)^{k+1}_{\mbfu,t})$, where we choose~$\R^-_0$ as branch cut for the complex logarithm in the variables~$\mbfu$ and~$t$. The same holds for~$\mcO_{\Lambda_\boldalpha}$ if~$\length{\boldalpha}\not\in \Z$.

\begin{lem}[Raising/lowering actions]\label{l:Calcoli}
The operators~\eqref{eq:DiffOpPhi} satisfy, for~$i,j\in[k]$, $j\neq i$,
\begin{equation}\label{eq:t:dsa1}
\begin{aligned}
E_{\alpha_i} f_{\boldalpha}=& \alpha_i f_{\boldalpha+\mbfe_i} \comm & E_{-\alpha_i} f_{\boldalpha}=&(1-\length{\boldalpha})f_{\boldalpha-\mbfe_i} \comm\\
E_{\alpha_i,-\alpha_j} f_{\boldalpha}=& \alpha_i f_{\boldalpha+\mbfe_i-\mbfe_j} \comm & J_{\alpha_i} f_{\boldalpha}=&(\length{\boldalpha}+\alpha_i-1) f_{\boldalpha} \fstop
\end{aligned}
\end{equation}

\begin{proof} The statement on~$J_{\alpha_i}$ is straightforward. Moreover,
\begin{align*}
E_{\alpha_i,-\alpha_j} f_{\boldalpha}=& \mbfu^\boldalpha t^{\length{\boldalpha}} \tonde{\sum_{\mbfm\geq \zero} \frac{\Poch{\boldalpha}{\mbfm} (m_i+\alpha_i) \mbfs^\mbfm}{\Poch{\length{\boldalpha}}{\length{\mbfm}} \mbfm!} - 
\frac{\Poch{\boldalpha}{\mbfm} m_i \mbfs^{\mbfm-\mbfe_i+\mbfe_j}}{\Poch{\length{\boldalpha}}{\length{\mbfm}} \mbfm!} }\\
=& \mbfu^\boldalpha t^{\length{\boldalpha}} \tonde{\sum_{\mbfm\geq \zero} \frac{\Poch{\boldalpha}{\mbfm} (m_i+\alpha_i) \mbfs^\mbfm}{\Poch{\length{\boldalpha}}{\length{\mbfm}} \mbfm!} - 
\frac{\Poch{\boldalpha}{\mbfm+\mbfe_i-\mbfe_j} (m_i+1)\mbfs^{\mbfm}}{\Poch{\length{\boldalpha}}{\length{\mbfm}} (\mbfm+\mbfe_i-\mbfe_j)!} }\\
=& \mbfu^\boldalpha t^{\length{\boldalpha}} \frac{\alpha_i}{\alpha_j-1} 
\tonde{\sum_{\mbfm\geq \zero} \frac{\Poch{\boldalpha+\mbfe_i-\mbfe_j}{\mbfm-\mbfe_i+\mbfe_j} (m_i+\alpha_i) \mbfs^\mbfm}{\Poch{\length{\boldalpha}}{\length{\mbfm}} \mbfm!} - 
\frac{\Poch{\boldalpha+\mbfe_i-\mbfe_j}{\mbfm} \mbfs^{\mbfm}}{\Poch{\length{\boldalpha}}{\length{\mbfm}} (\mbfm-\mbfe_j)!} }\\
=& \mbfu^\boldalpha t^{\length{\boldalpha}} \frac{\alpha_i}{\alpha_j-1}
\tonde{\sum_{\mbfm\geq \zero} \frac{\Poch{\boldalpha+\mbfe_i-\mbfe_j}{\mbfm} (m_j+\alpha_j-1) \mbfs^\mbfm}{\Poch{\length{\boldalpha}}{\length{\mbfm}} \mbfm!} - 
\frac{\Poch{\boldalpha+\mbfe_i-\mbfe_j}{\mbfm} m_j \mbfs^{\mbfm}}{\Poch{\length{\boldalpha}}{\length{\mbfm}} \mbfm!} }\\
=&\alpha_i f_{\boldalpha+\mbfe_i-\mbfe_j} \comm
\\
%
E_{\alpha_i} f_{\boldalpha}=& \mbfu^\boldalpha t^{\length{\boldalpha}}\tonde{\sum_{\mbfm\geq \zero} \frac{\Poch{\boldalpha}{\mbfm} (m_i+\alpha_i) \mbfs^\mbfm}{\Poch{\length{\boldalpha}}{\length{\mbfm}} \mbfm!} - 
\frac{\Poch{\boldalpha}{\mbfm} m_i(\length{\mbfm}-1) \mbfs^{\mbfm-\mbfe_i}}{\Poch{\length{\boldalpha}}{\length{\mbfm}} \mbfm!} }\\
=&\mbfu^\boldalpha t^{\length{\boldalpha}}\tonde{\sum_{\mbfm\geq \zero} \frac{\Poch{\boldalpha}{\mbfm+\mbfe_i} \mbfs^\mbfm}{\Poch{\length{\boldalpha}}{\length{\mbfm}} \mbfm!} - 
\frac{\Poch{\boldalpha}{\mbfm} (\length{\mbfm}-1) \mbfs^{\mbfm-\mbfe_i}}{\Poch{\length{\boldalpha}}{\length{\mbfm}} (\mbfm-\mbfe_i)!} }\\
=&\mbfu^\boldalpha t^{\length{\boldalpha}}\tonde{\sum_{\mbfm\geq \zero} \frac{\Poch{\boldalpha}{\mbfm+\mbfe_i} \mbfs^\mbfm}{\Poch{\length{\boldalpha}}{\length{\mbfm}} \mbfm!} - 
\frac{\Poch{\boldalpha}{\mbfm+\mbfe_i} \length{\mbfm} \mbfs^{\mbfm}}{\Poch{\length{\boldalpha}}{\length{\mbfm}+1} \mbfm!} }\\
=&\mbfu^\boldalpha t^{\length{\boldalpha}}\frac{\alpha_i}{\length{\boldalpha}}\tonde{\sum_{\mbfm\geq \zero} \frac{\Poch{\boldalpha+\mbfe_i}{\mbfm} \mbfs^\mbfm}{\Poch{\length{\boldalpha}+1}{\length{\mbfm}-1} \mbfm!} - 
\frac{\Poch{\boldalpha+\mbfe_i}{\mbfm} \length{\mbfm} \mbfs^{\mbfm}}{\Poch{\length{\boldalpha}+1}{\length{\mbfm}} \mbfm!} }\\
=&\mbfu^\boldalpha t^{\length{\boldalpha}}\frac{\alpha_i}{\length{\boldalpha}}\tonde{\sum_{\mbfm\geq \zero} \frac{\Poch{\boldalpha+\mbfe_i}{\mbfm} \mbfs^\mbfm (\length{\boldalpha}+\length{\mbfm})}{\Poch{\length{\boldalpha}+1}{\length{\mbfm}} \mbfm!} - 
\frac{\Poch{\boldalpha+\mbfe_i}{\mbfm} \length{\mbfm} \mbfs^{\mbfm}}{\Poch{\length{\boldalpha}+1}{\length{\mbfm}} \mbfm!} }\\
=&\alpha_i f_{\boldalpha+\mbfe_i} \comm
\\
%
E_{-\alpha_i} f_{\boldalpha}=& \mbfu^\boldalpha t^{\length{\boldalpha}}\tonde{\sum_{\mbfm\geq \zero} \frac{\Poch{\boldalpha}{\mbfm} \mbfs^{\mbfm+\mbfe_i}}{\Poch{\length{\boldalpha}}{\length{\mbfm}} \mbfm!} -
\frac{\Poch{\boldalpha}{\mbfm} \mbfs^{\mbfm}}{\Poch{\length{\boldalpha}}{\length{\mbfm}} \mbfm!}(\length{\mbfm}+\length{\boldalpha}-1) }\\
=& \mbfu^\boldalpha t^\length{\boldalpha} \tonde{\sum_{\mbfm\geq \zero} \frac{\Poch{\boldalpha}{\mbfm-\mbfe_i} m_i \mbfs^\mbfm}{\Poch{\length{\boldalpha}}{\length{\mbfm}-1} \mbfm!} -
\frac{\Poch{\boldalpha}{\mbfm} \mbfs^{\mbfm}}{\Poch{\length{\boldalpha}}{\length{\mbfm}} \mbfm!} (\length{\mbfm}+\length{\boldalpha}-1) }\\
=& \mbfu^\boldalpha t^\length{\boldalpha} \frac{\length{\boldalpha}-1}{\alpha_i-1}\times\\
&\times
\tonde{\sum_{\mbfm\geq \zero} \frac{\Poch{\boldalpha-\mbfe_i}{\mbfm} m_i \mbfs^\mbfm}{\Poch{\length{\boldalpha}-1}{\length{\mbfm}} \mbfm!} -
\frac{\Poch{\boldalpha-\mbfe_i}{\mbfm} (\alpha_i+m_i-1) \mbfs^{\mbfm}}{\Poch{\length{\boldalpha}-1}{\length{\mbfm}} (\length{\boldalpha}+\length{\mbfm}-1) \mbfm!} (\length{\mbfm}+\length{\boldalpha}-1) }\\
=&(1-\length{\boldalpha})f_{\boldalpha-\mbfe_i} \fstop \qedhere
\end{align*}
\end{proof}
\end{lem}

\begin{rem}\label{r:OpDep}
The variables~$\mbfu$ and $t$ are merely auxiliary (cf.~\cite[\S1]{Mil73b}).
The operators do not depend on the parameter~$\boldalpha$, rather, the subscripts indicate which indices they affect.
Heuristically, the action of the operators~\eqref{eq:DiffOpPhi} given in Lemma~\ref{l:Calcoli} may be derived from that~\cite[(1.5)]{Mil72} of operators in the dynamical symmetry algebra of~${}_kF_D$ by a formal contraction procedure~\cite[p.~1398]{Mil72}, letting (in the notation of~\cite{Mil72})~$\alpha=0$, $\boldbeta=\boldalpha$, $\gamma=\length{\boldalpha}$ and dropping redundancies.
\end{rem}

\begin{rem}
If~$\length{\boldalpha}=1$, the action of the lowering operators~$E_{-\alpha_i}$ vanishes. This is natural when regarding~$f_\boldalpha$ as a formal power series, whereas it is conventional when regarding~$f_\boldalpha$ as a meromorphic function, for the functions~$(1-\length{\boldalpha})f_{\boldalpha-\mbfe_i}$ are in fact ---~after cancellations~--- well-defined, not identically vanishing, and holomorphic in~$\mbfs$ even for~$\length{\boldalpha}=1$. The convention here reads~$0\times\infty=0$, which is consistent with the usual convention in measure theory when we identify~$\length{\boldalpha}-1$ with the quantity~$(\sigma-\delta_y)X$ for any~$y$ in~$X$; the reason for such identification will be apparent in~\S\ref{ss:InfDimTDue} below.
\end{rem}

\begin{cor}\label{c:Fix}
The operators~\eqref{eq:DiffOpPhi} fix~$\mcO_{\Lambda_\boldalpha}$ for any~$\boldalpha\in \C^k$.
\end{cor}

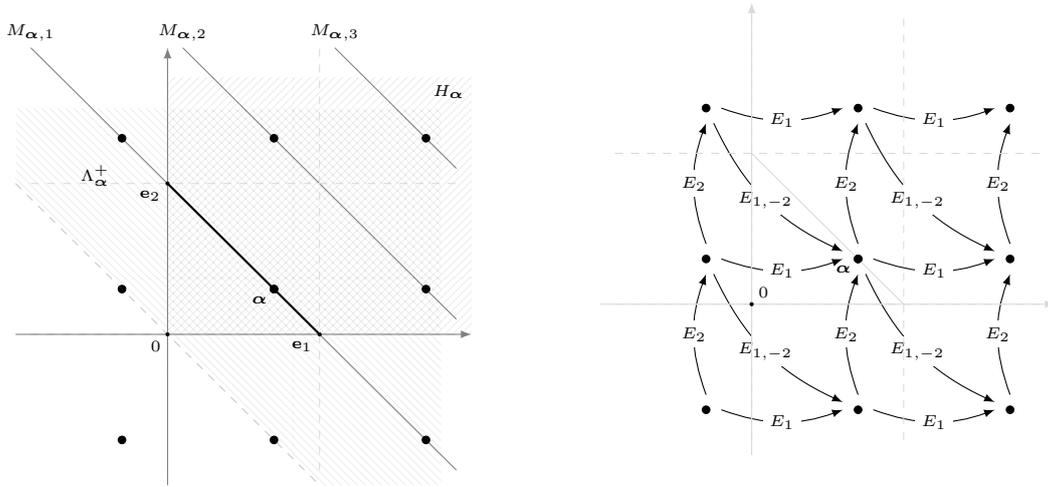
\begin{figure}[ht]
    \centering
    \begin{subfigure}[b]{0.49\textwidth}
    \centering
       \begin{tikzpicture}[scale=2]
\tiny

\def\ax{.7} \def\ay{.3}
\def\inset{.1}
\def\boxes{1}
\def\rboxes{1}

\draw[color=gray] (-\boxes+\inset,\boxes+1-\inset) node[black, above]{$M_{\boldalpha,1}$} -- (\boxes+1-\inset,-\boxes+\inset);
\draw[color=gray] (-\boxes+1+\inset,\boxes+1-\inset) node[black, above]{$M_{\boldalpha,2}$} -- (\boxes+1-\inset,-\boxes+\inset+1);
\draw[color=gray] (-\boxes+2+\inset,\boxes+1-\inset) node[black, above]{$M_{\boldalpha,3}$} -- (\boxes+1-\inset,-\boxes+\inset+2);

\foreach \x in {-\boxes,...,\boxes}
    \foreach \y in {-\boxes,...,\boxes}
	\draw[fill] (\ax+\x,\ay+\y) circle (.025);
%
%
%
		
\draw[thick] (0,1) -- (1,0);

\draw[draw=none, fill, pattern=north east lines, pattern color=grey1, opacity=.3] (0,\boxes+1-3*\inset) -- (0,0) -- (\boxes+1, 0) -- (\boxes+1, \boxes+1-3*\inset) -- cycle;

\draw[draw=none, fill, pattern=north west lines, pattern color=grey1, opacity=.3] (-\boxes
,\boxes) -- (-\boxes
,\boxes+1-5*\inset) -- (\boxes+1-2*\inset, \boxes+1-5*\inset) -- (\boxes+1-2*\inset, -\boxes
)
 -- (\boxes, -\boxes
 ) -- cycle;

\draw[dashed, color=grey1, opacity=.3] (-\boxes,\boxes) -- (\boxes, -\boxes);

\draw (-\boxes+\inset/2,\boxes+1-5*\inset)+(1/3,-1/3) node[below right] {$\Lambda^+_\boldalpha$} ;
\draw (\boxes+1, \boxes+1-3*\inset) node[below left] {$H_\boldalpha$} ;

\draw[help lines, color=gray!30, dashed] (-\boxes+\inset,-\boxes+\inset) grid (\boxes-\inset+1,\boxes-\inset+1);
\draw[->, color=gray] (-\boxes,0)--(\boxes+1,0); 
\draw[->, color=gray] (0,-\boxes)--(0,\boxes+1-\inset); 

\draw[fill] (\ax,\ay) node[below left]{$\boldalpha$} circle (.025);
\draw[fill] (0,0) node[below left]{$0$} circle (.01);
\draw[fill] (0,1) node[below left]{$\mbfe_2$} circle (.01);
\draw[fill] (1,0) node[below left]{$\mbfe_1$} circle (.01);
\end{tikzpicture}
    \end{subfigure}%
    \begin{subfigure}[b]{0.49\textwidth}
    \centering
        \raisebox{.4cm}
{\begin{tikzpicture}[scale=2]
\tiny

\def\ax{.7} \def\ay{.3}
\def\inset{.1}
\def\boxes{1}
\def\rboxes{1}


\foreach \x in {-\boxes,...,\boxes}
    \foreach \y in {-\boxes,...,\boxes}
	\draw[fill] (\ax+\x,\ay+\y) circle (.025);
%
\foreach \x in {-\rboxes,...,\number\numexpr\rboxes-1}
	\foreach \y in {-\rboxes,...,\rboxes}
		\draw[->] (\ax+\x+\inset,\ay+\y) to [bend right=20] node[midway, fill=white] {$E_1$} (\ax+\x+1-\inset, \ay+\y);

\foreach \x in {-\rboxes,...,\number\numexpr\rboxes}
	\foreach \y in {-\rboxes,...,\number\numexpr\rboxes-1}
		\draw[->] (\ax+\x,\ay+\y+\inset) to [bend left=20] node[midway,fill=white] {$E_2$} (\ax+\x, \ay+\y+1-\inset);
				
\foreach \x in {0,...,\number\numexpr\rboxes}
	\foreach \y in {-\rboxes,...,\number\numexpr\rboxes-1}
		\draw[->] (\ax+\x-1+\inset/2,\ay+\y+1-\inset) to [bend right=20] node[midway, fill=white] {$\displaystyle{E_{1,-2}}$} (\ax+\x-\inset, \ay+\y+\inset/2);
		
\draw[color=gray!30] (0,1) -- (1,0);

\draw[help lines, color=gray!30, dashed] (-\boxes+\inset,-\boxes+\inset) grid (\boxes-\inset+1,\boxes-\inset+1);
\draw[->, color=gray!30] (-\boxes,0)--(\boxes+1,0); 
\draw[->, color=gray!30] (0,-\boxes)--(0,\boxes+1); 

\draw[fill] (\ax,\ay) node[below left]{$\boldalpha$} circle (.025);
\draw[fill] (0,0) node[above right]{$0$} circle (.01);
\end{tikzpicture}
}
    \end{subfigure}
    \caption{
    \emph{(both)} Each marked point corresponds to some~$\boldalpha'\in\Lambda_\boldalpha$ for fixed~$\boldalpha$, and is chosen to indicate the one-dimensional vector space~$\mcO_{\boldalpha'}$.
    \emph{(left)} The gray anti-diagonal lines denote the \emph{isoplethic surfaces}: marked points~$\boldalpha'$ lying on these surfaces belong to~$M_{\boldalpha,\ell}$, i.e. they have fixed length~$\length{\boldalpha'}=\ell\in\N_1$. The simplex~$\Delta^1$ is marked as a thick black segment. Analogously, marked points lying in the North-West dashed region delimited by the hyper-plane of equation~$\length{\mbfy}=0$ belong to the semi-lattice~$\Lambda_\boldalpha^+$, whereas marked points lying in the first hyper-octant (in the figure: the North-East dashed quadrant) belong to~$H_\boldalpha$.
    \emph{(right)} The action of operators in~$\rho_\boldalpha(\mfn^+_2)$ on the lattice~$\mcO_{\Lambda_\boldalpha}$ for~$\boldalpha=\seq{\tfrac 23,\tfrac13}$ is shown.
}
\label{fig:LMH}
\end{figure}

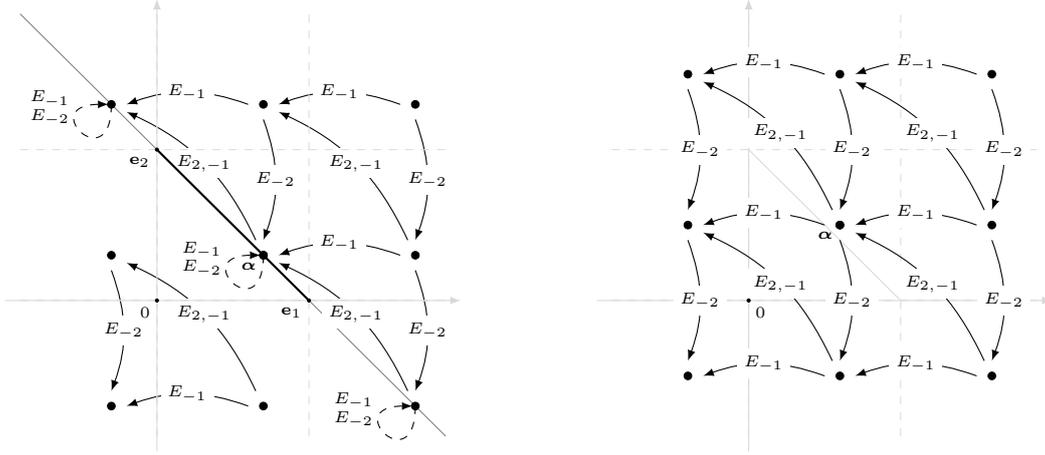
\begin{figure}[ht]
    \centering
    \begin{subfigure}[b]{0.49\textwidth}
    \centering
       \begin{tikzpicture}[scale=2]
\tiny

\def\ax{.7} \def\ay{.3}
\def\inset{.1}
\def\boxes{1}
\def\rboxes{1}

\draw[color=gray] (-\boxes+\inset,\boxes+1-\inset) -- (\boxes+1-\inset,-\boxes+\inset);

\foreach \x in {-\boxes,...,\boxes}
    \foreach \y in {-\boxes,...,\boxes}
	\draw[fill] (\ax+\x,\ay+\y) circle (.025);
%
\foreach \x in {-\rboxes,...,\number\numexpr\rboxes-1}
	\foreach \y in {-\x,...,\rboxes}
		\draw[->] (\ax+\x+1-\inset,\ay+\y) to [bend right=20] node[midway, fill=white] {$E_{-1}$} (\ax+\x+\inset, \ay+\y);

\foreach \x in {0,...,\number\numexpr\rboxes-1}
	\foreach \y in {-\rboxes,...,\number\numexpr-\x-1}
		\draw[->] (\ax+\x-\inset,\ay+\y) to [bend right=20] node[midway, fill=white] {$E_{-1}$} (\ax+\x-1+\inset, \ay+\y);

\foreach \x in {-\rboxes,...,\number\numexpr\rboxes-1}
	\foreach \y in {-\x,...,\rboxes}
		\draw[->] (\ax+\x+1,\ay+\y-\inset) to [bend left=20] node[midway,fill=white] {$E_{-2}$} (\ax+\x+1, \ay+\y-1+\inset);

\foreach \x in {-\rboxes,...,\number\numexpr\rboxes-1-1}
	\foreach \y in {-\rboxes,...,\number\numexpr-\x-1-1}
		\draw[->] (\ax+\x,\ay+\y+1-\inset) to [bend left=20] node[midway,fill=white] {$E_{-2}$} (\ax+\x, \ay+\y+\inset);
		
\foreach \x in {-\rboxes,...,\number\numexpr\rboxes-1}
	\foreach \y in {-\x,...,\rboxes}
		\draw[->] (\ax+\x+1-\inset/2,\ay+\y-1+\inset) to [bend right=20] node[midway, fill=white] {$\displaystyle{E_{2,-1}}$} (\ax+\x+\inset, \ay+\y-\inset/2);

\foreach \x in {0,...,\number\numexpr\rboxes-1}
	\foreach \y in {-\rboxes,...,\number\numexpr-\x-1}
		\draw[->] (\ax+\x-\inset/2,\ay+\y+\inset) to [bend right=20] node[midway, fill=white] {$\displaystyle{E_{2,-1}}$} (\ax+\x-1+\inset, \ay+\y+1-\inset/2);
		
\foreach \x in {-\rboxes,...,\rboxes}
	\path[dashed, ->] (\ax-\x, \ay+\x) edge [in=180,out=-90, loop] node[above left, text width=.6cm] {$E_{-1}$ $E_{-2}$} (\ax-\x, \ay+\x);

\draw[thick] (0,1) -- (1,0);

\draw[help lines, color=gray!30, dashed] (-\boxes+\inset,-\boxes+\inset) grid (\boxes-\inset+1,\boxes-\inset+1);
\draw[->, color=gray!30] (-\boxes,0)--(\boxes+1,0); 
\draw[->, color=gray!30] (0,-\boxes)--(0,\boxes+1); 

\draw[fill] (\ax,\ay) node[below left]{$\boldalpha$} circle (.025);
\draw[fill] (0,0) node[below left]{$0$} circle (.01);
\draw[fill] (0,1) node[below left]{$\mbfe_2$} circle (.01);
\draw[fill] (1,0) node[below left]{$\mbfe_1$} circle (.01);
\end{tikzpicture}
    \end{subfigure}%
    \begin{subfigure}[b]{0.49\textwidth}
    \centering
        \raisebox{.4cm}
{\begin{tikzpicture}[scale=2]
\tiny

\def\ax{3/5} \def\ay{1/2}
\def\inset{.1}
\def\boxes{1}
\def\rboxes{1}

\draw[color=gray!30] (0,1) -- (1,0);

\draw[help lines, color=gray!30, dashed] (-\boxes+\inset,-\boxes+\inset) grid (\boxes-\inset+1,\boxes-\inset+1);
\draw[->, color=gray!30] (-\boxes,0)--(\boxes+1,0); 
\draw[->, color=gray!30] (0,-\boxes)--(0,\boxes+1); 


\foreach \x in {-\boxes,...,\boxes}
    \foreach \y in {-\boxes,...,\boxes}
	\draw[fill] (\ax+\x,\ay+\y) circle (.025);
%
\foreach \x in {-\rboxes,...,\number\numexpr\rboxes-1}
	\foreach \y in {-\rboxes,...,\rboxes}
		\draw[->] (\ax+\x+1-\inset,\ay+\y) to [bend right=20] node[midway, fill=white] {$E_{-1}$} (\ax+\x+\inset, \ay+\y);

\foreach \x in {-\rboxes,...,\number\numexpr\rboxes}
	\foreach \y in {-\rboxes,...,\number\numexpr\rboxes-1}
		\draw[->] (\ax+\x,\ay+\y+1-\inset) to [bend left=20] node[midway,fill=white] {$E_{-2}$} (\ax+\x, \ay+\y+\inset);
		
\foreach \x in {-\rboxes,...,\number\numexpr\rboxes-1}
	\foreach \y in {-\rboxes,...,\number\numexpr\rboxes-1}
		\draw[->] (\ax+\x+1-\inset/2,\ay+\y+\inset) to [bend right=20] node[midway, fill=white] {$\displaystyle{E_{2,-1}}$} (\ax+\x+\inset, \ay+\y+1-\inset/2);

\draw[fill] (\ax,\ay) node[below left]{$\boldalpha$} circle (.025);
\draw[fill] (0,0) node[below right]{$0$} circle (.01);
\end{tikzpicture}
}
    \end{subfigure}
    \caption{
    \emph{(both)} Each marked point corresponds to some~$\boldalpha'\in\Lambda_\boldalpha$ and is chosen to indicate the one-dimensional vector space~$\mcO_{\boldalpha'}$.
    \emph{(left)} The action of operators in~$\rho_\boldalpha(\mfn^-_2)$ on the lattice~$\mcO_{\Lambda_\boldalpha}$ for~$\boldalpha=\seq{\tfrac 23,\tfrac13}$ is shown. Since~$\length{\boldalpha}\in \Z$, the lowering operators~$E_{-\alpha_1}, E_{-\alpha_2}$ \emph{(left)} vanish identically on the lowest positive isoplethic line (in gray), containing the standard simplex (the thick segment): their action is here represented by a dashed loop.
    \emph{(right)} The action of operators in~$\rho_\boldalpha(\mfn^-_2)$ on the lattice~$\mcO_{\Lambda_\boldalpha}$ for~$\boldalpha=\seq{\tfrac 35,\tfrac12}$ is shown. Since~$\length{\boldalpha}\in \R\setminus \Z$, the lowering operators~$E_{-\alpha_1}, E_{-\alpha_2}$ never vanish.
    }
\label{fig:Lowering}
\end{figure}

In the statement of the next Lemma and in the diagrams in Fig.s~\ref{fig:LMH} and~\ref{fig:Lowering} we write for simplicity~$E_i$ in place of~$E_{\alpha_i}$ and analogously for all other operators.
\begin{lem}\label{l:Commu}
For~$\boldalpha\in\C^k$ consider the operators in~$\mfg_k$ as restricted to~$\mcO_{\Lambda_\boldalpha}$. The following commutation relations hold:
\begin{align*}
[J_{i}-J_{j},E_{p,-q}]=&\begin{cases}
2E_{p,-q} &\textrm{if}\;\; i=p,j=q\\
-2E_{p,-q} &\textrm{if}\;\; i=q,j=p\\
E_{p,-q} &\textrm{if}\;\; i=p,j\neq q \\
&\;\;\textrm{or}\;\; i\neq p,j=q\\
-E_{p,-q} &\textrm{if}\;\; i=q,j\neq p \\
&\;\;\textrm{or}\;\; i\neq q,j= p\\
\zero &\textrm{otherwise}
\end{cases}\;,&
\begin{aligned}
[E_{i,-j},E_{p,-q}]=&\begin{cases} 
J_{i}-J_{j} & \textrm{if}\;\; i=q,j=p\\
-E_{p,-j} & \textrm{if}\;\; i=q, j\neq p\\
E_{i,-q} & \textrm{if}\;\; i\neq q, j=p\\
\zero & \textrm{otherwise}
\end{cases}
\\
[E_{i},E_{-p}]=&
\begin{cases}
J_{i} &\textrm{if}\;\; i=p\\
E_{i,-p} &\textrm{if}\;\; i\neq p
\end{cases}
\\
[J_{i},E_{\pm p}]=&\pm2\delta_{ip}E_{\pm p} \comm
\end{aligned}\;,
\end{align*}
\begin{align*}
[J_{i}-J_{j}, J_{p}]=[E_{i},E_{-p}]=\zero \comm
\end{align*}
where~$i,j=0,\dotsc, k$ and~$p,q=1,\dotsc,k$ with~$i\neq j$, $p\neq q$ and, conventionally,
\begin{align*}
J_{0}=\zero\comm E_{0,-j}=E_{-j}\comm E_{i,0}=E_{i}\comm 0\neq 0 \fstop
\end{align*}

\begin{proof} Given the action of the operators in~\eqref{eq:t:dsa1} straightforward computations yield 
\begin{align*}
[J_{\alpha_i}-J_{\alpha_j},E_{\alpha_p,-\alpha_q}]f_\boldalpha=&\tonde{\delta_{ip}-\delta_{iq}-\delta_{jp}+\delta_{jq}} \alpha_p f_{\boldalpha+\mbfe_p-\mbfe_q}\comm
\\
[E_{\alpha_i,-\alpha_j},E_{\alpha_p,-\alpha_q}]f_\boldalpha=&\tonde{(\delta_{ip}-\delta_{iq})\alpha_p-(\delta_{ip}-\delta_{jp})\alpha_i}f_{\boldalpha+\mbfe_i-\mbfe_j+\mbfe_p-\mbfe_q}\comm
\\
[E_{\alpha_i,E_{-\alpha_p}}]f_\boldalpha=&(\alpha_i+(\length{\boldalpha}-1)\delta_{ip})f_{\boldalpha+\mbfe_i-\mbfe_p}\comm
\\
[J_{\alpha_i},E_{\alpha_p}]f_{\boldalpha}=&2\delta_{ip}\alpha_pf_{\boldalpha+\mbfe_p}\comm
\\
[J_{\alpha_i},E_{-\alpha_p}]f_{\boldalpha}=&-2\delta_{ip}(1-\length{\boldalpha})f_{\boldalpha-\mbfe_p} \fstop \qedhere
\end{align*}
\end{proof}
\end{lem}

\begin{prop}\label{p:Repr}
Let~$\rho\colon \mfl_k\rar \End (\mcO)$ be the linear map defined by
\begin{align*}
e_{0,i}\mapsto E_{\alpha_i}\comm e_{i,j}\mapsto E_{\alpha_i,-\alpha_j}\comm h_{0,i}\mapsto J_{\alpha_i}\comm f_{j,0}\mapsto E_{-\alpha_j}\comm f_{j,i}\mapsto E_{\alpha_j,-\alpha_i}
\end{align*}
where~$i,j\in [k]$ with~$j>i$.
Then, for any fixed~$\boldalpha\in \C^k$, the pair~$\rho_\boldalpha\eqdef (\rho(\emparg)\restr_{\mcO_{\Lambda_\boldalpha}}, \mcO_{\Lambda_\boldalpha})$ is a faithful Lie algebra representation of~$\mfl_k$ with image~$\mfg_k\restr_{\mcO_{\Lambda_\boldalpha}}$. Furthermore, the functions~$f_\boldalpha$ transform as basis vectors for~$\rho_\boldalpha$, in the sense that for every~$v$ in the basis for~$\mfl_k$ and every~$\boldalpha'$ in~$\Lambda_\boldalpha$ there exists a unique~\mbox{$\boldalpha''=\boldalpha''(\boldalpha', v)$} in~$\Lambda_\boldalpha$ such that $(\rho_\boldalpha v)\, f_{\boldalpha'} \propto f_{\boldalpha''}$.

\begin{proof}
By Corollary~\ref{c:Fix},~$\rho_\boldalpha$ is a well-defined linear morphism into~$\End(\mcO_{\Lambda_\boldalpha})$. The fact that~$f_\boldalpha$ transforms as a basis vector of~$\mcO_{\Lambda_\boldalpha}$ is an immediate consequence of Lemma~\ref{l:Calcoli}. For~\mbox{$\boldalpha'\in \Lambda_\boldalpha$} such that~$\Re^\compo\boldalpha'>\uno$, the actions of operators in~\eqref{eq:DiffOpPhi} on~$\mcO_{\boldalpha'}$ are mutually different again by Lemma~\ref{l:Calcoli}, hence~$\rho_\boldalpha$ is injective. In order to show that~$\rho_\boldalpha \mfl=\mfg\restr_{\mcO_{\Lambda_\boldalpha}}$ is a Lie algebra of type~$A_k$ and that~$\rho_\boldalpha$ is a Lie algebra representation, it suffices to verify Serre relations~\cite[\S18.1]{Hum72} of type~$A$ for the operators~$\rho_\boldalpha v$ with~$v=v_j$ in an $\mfsl_2$-triple corresponding to the simple root~$\gamma_j$ in~$\Psi_k$. These are readily deduced from Lemma~\ref{l:Commu}.
\end{proof}
\end{prop}

\begin{thm}\label{c:Due}
For~$\boldalpha$ in~$\interior\Delta^{k-1}$ and~$\mbfp\in (\Z_0^+)^k$ denote by~$\Dir{\boldalpha}^\mbfp$ the posterior distribution of~$\Dir{\boldalpha}$ given atoms of mass~$p_i$ at point~$i\in [k]$. 
Then,
\begin{enumerate}[$(i)$]
\item\label{i:c:Due1} the semi-lattice~$\mcO_{\Lambda^+_\boldalpha}$ is a weight $\mfl$-module and~$\mfU(\mfl)$-module;

\item\label{i:c:Due2} the action of the universal enveloping algebra~$\mfU(\mff)<\mfU(\mfl)$ globally fixes~$\mcO_{M_{\boldalpha,\ell}}$ for all~$\ell\in \Z^+$, while the action of~$\mfU(\mfh)\oplus\mfU(\mfr^+)$ globally fixes~$\mcO_{H_{\boldalpha}}$;

\item\label{i:c:Due3} for every~$\mbfp\in (\Z^+_0)^k$ there exists a unique~$v= v(\mbfp)\in \mfU(\mfr^+)$ such that~$v.\mcO_\boldalpha\cong\C \FT{\Dir{\boldalpha}^\mbfp}$;

\item\label{i:c:Due4} the canonical action of~$\mfS_k$ on~$\msP([k])$ corresponds to the natural action of the unique subgroup (isomorphic to)~$\mfS_k$ of the Weyl group of~$\mfl_k$ permuting roots corresponding to basis elements in~$\mfr^+_k$. 
\end{enumerate}

\begin{proof}
By~\eqref{eq:t:dsa1}, the operators~$\rho\, \mfsl_{k+1}(\C)$ fix~$\mcO_{\Lambda_\boldalpha^+}\subset \mcO$, thus~$\rho_\boldalpha$ is a (faithful) Lie algebra representation by Proposition~\ref{p:Repr}, hence~$\mcO_{\Lambda_\boldalpha^+}\subset \mcO$ is an~$\mfl$-module for the  linear extension of the action~$v.f_{\boldalpha'}\eqdef (\rho_\boldalpha v) f_{\boldalpha'} $ varying~$v$ in the basis of~$\mfl$. The extension to a representation of~$\mfU(\mfl)$ is standard from the universal property of universal enveloping algebras (e.g.~\cite[\S17.2]{Hum72}).

In order to prove~\iref{i:c:Due1}-\iref{i:c:Due2} it suffices to show that, for all~$\boldalpha'\in \Lambda^+_\boldalpha$ and~$\ell\in\Z^+$, one has
\begin{align*}
h_{0,i}.f_{\boldalpha'}=(\length{\boldalpha'}-1+\alpha_i')f_{\boldalpha'}\comm v.\mcO_{M_{\boldalpha',\ell}}\subset \mcO_{M_{\boldalpha',\ell}}\comm w. \mcO_{H_\boldalpha}\subset \mcO_{H_\boldalpha}
\end{align*}
for all $i\in [k]$,~$v$ in the basis of~$\mff$,~$\ell\in \N_1$ and~$w$ in the basis for~$\mfh\oplus \mfr^+$. All of the above follow immediately from Lemma~\ref{l:Calcoli}. Notably, since~$\length{\boldalpha}=1$,~$\mfh$ acts on~$\mcO_\boldalpha$ precisely by weight~$\boldalpha$.

Since~$\boldalpha\in \Delta^{k-1}$, then~$f_{\boldalpha+\mbfp}(\emparg,\uno,1)=\FT{\Dir{\boldalpha+\mbfp}}(\emparg)$. By the Bayesian property of~$\Dir{\boldalpha}$ the space~$\mcO_{H_\boldalpha}$ is spanned precisely by the Fourier transforms of the form~$\FT{\Dir{\boldalpha}^\mbfp}$. It remains to show that~$\mfU(\mfr^+).\mcO_\boldalpha=\mcO_{H_\boldalpha}$.
Setting~$v=e_1^{p_1}\cdots e_k^{p_k}\in \mfU(\mfr^+_k)$ yields~$v.\mcO_{\boldalpha}=\mcO_{\boldalpha+\mbfp}$ as required.
The uniqueness of~$v$ follows by the fact that, since~$\mfr^+$ is Abelian,~$\mfU(\mfr^+)$ coincides with the (Abelian) symmetric algebra generated by~$\mfr^+$ (see~\cite[\S17.2]{Hum72}). This proves~\iref{i:c:Due3}.

\smallskip

In order to show~\iref{i:c:Due4}, recall (e.g.~\cite[\S12.1]{Hum72}) that the Weyl group~$W_k$ of~$\Psi_k$ is isomorphic to~$\mfS_{k+1}$ and its action on~$\Psi_k$ may be canonically identified as dual to the action of~$\mfS_{k+1}$ on~$\mfh_k$ via conjugation by permutation matrices in~$\mfP_{k+1}\cong\mfS_{k+1}<GL(\mfh_k)\cong GL_{k+1}(\C)$.
Let~$\mfP_{2:k+1}<GL_{k+1}(\C)$ denote the subgroup of permutations matrices whose action on~$\mathrm{Mat}_{k+1}(\C)$ fixes the first row and column. Clearly~$\mfS_k\cong\mfP_{2:k+1}<\mfP_{k+1}$. Composing the isomorphism~$\rho_\boldalpha$ with the identification of the action of~$\mfP_{k+1}$ above completes the proof.
\end{proof}
\end{thm}

\subsection{Invariant measures on simplices and affine spheres}\label{ss:Algebra}
In the following we lay out a comparison between the results obtained in the previous section and some known facts about the \emph{multiplicative infinite-dimensional Lebesgue measure $\mcL^+$}~\cite{TsiVerYor00}. For the reader's convenience, let us briefly recall the construction of~$\mcL^+$ given in~\cite{Ver07}.

\paragraph{The measure~$\mcL^+_{\beta,\sigma}$}
For~$\mbfy\eqdef\seq{y_1,\dotsc,y_k}\in\R^k$ let~$\diag(\mbfy)$ be the diagonal matrix with entries~$\mbfy$ and set
\begin{align*}
M^{k-1}_r\eqdef \set{\mbfy\in\R^k_+\mid \mbfy^{\uno}=r} \qquad r>0\comm
\end{align*}
hereby termed the $(k-1)$-dimensional \emph{affine sphere} of radius $r$. Set~$dSL^+_k(\R)\eqdef\diag(M^{k-1}_1)$, namely, the connected component of the identity in the maximal toral subgroup (of positive definite diagonal matrices) in the special linear group~$SL_k(\R)$. Let further~$\ell^k_\emparg$ denote the natural self-action of~$dSL^+_k(\R)$ by left-multiplication and notice that it induces an action~$\eqdef dSL^+_k(\R)\acts M^{k-1}_r$ given by~$\diag^{-1} \circ \ell^k_\emparg\circ \diag$, in the following also denoted by~$\ell^k_\emparg$. Since~$\ell^k_\emparg$ is free and transitive,~$M^{k-1}_r$ admits an $\ell^k_\emparg$-invariant measure~$\lambda_{k,r}$.

Let now~$x_1,\dotsc,x_k$ be $X$-valued independent $\sigma$-distributed random variables, set~$\iota\colon \mbfy\mapsto \sum_i^k y_i\delta_{x_i}$ and let~$\mu_{k,r}\eqdef \iota_\pfwd \lambda_{k,r}$ be the image measure of~$\lambda_{k,r}$ induced on~$\Mbp(X)$ via~$\iota$.
In the case~$r_{k,\beta}\eqdef \exp(-\beta k^2)$, it was shown in~\cite[Thm.~3]{Ver07} that the sequence~$\seq{\mu_{k,r_{k,\beta}}}_k$ converges (in a suitable sense, see~\cite[\S3]{Ver07}) to the multiplicative infinite-dimensional Lebesgue measure~$\mcL^+_{\beta,\sigma}$, the unique measure on~$\Mbp(X)$ satisfying
\begin{align*}
\int_{\Mbp(X)}\diff \mcL^+_{\beta,\sigma}(\eta)\, \exp(\eta f)=\exp\tonde{-\beta\, \sigma\ln f} \qquad \ln f\in\mcC_c \fstop
\end{align*}

It was shown in~\cite{Ver07} that~$\mcL^+_{\beta,\sigma}\ll \GP_{\beta,\sigma}$ and that~$\mcL^+_{\beta,\sigma}$ is a positive $\sigma$-finite Borel measure on~$\Mbp(X)$.

\paragraph{The commutative action of~$dSL^+_k(\R)$} Let~$\mbfX\in\mfP_k(X,\tau(X),\sigma)$ and set~$\boldalpha\eqdef\beta\ev^\mbfX\sigma$. The finite-dimensional marginalizations~$L_\boldalpha\eqdef\ev^\mbfX_\pfwd \mcL^+_{\beta,\sigma}$ induced by~$\mbfX$ in the sense of~\eqref{eq:PartX} satisfy (see~\cite[Prop.~2]{Ver07})
\begin{align*}
L_\boldalpha=\car_{M^{k-1}_{r_{k,\beta}}}(\mbfy)\,\frac{\mbfy^{\boldalpha-\uno}}{\Gamma(\boldalpha)}\diff\mbfy \comm
\end{align*}
to be compared with the density~\eqref{eq:DensityD} of the Dirichlet distribution.
The $\ell^k_\emparg$-invariance of the measures~$\lambda_{k,r_{k,\beta}}$ on rescaled affine spheres corresponds $(a)$ in the finite-dimensional case~--- to the \emph{projective} invariance of the measures~$L_\boldalpha$ with respect to the same action, with Radon--Nikod\'ym derivative
\begin{align*}
\forallae{L_\boldalpha} \mbfy\in M_{r_{k,\beta}}^{k-1}\qquad\frac{\diff (\ell^k_\mbfs)_\pfwd L_\boldalpha}{\diff L_{\boldalpha}}(\mbfy)\equiv \mbfs^{-\boldalpha}\equiv \exp(-\boldalpha \cdot \ln^\compo\mbfs)\comm \diag(\mbfs)\in dSL^+_k(\R)\semicolon
\end{align*}
and $(b)$ in the infinite-dimensional case~--- to the \emph{projective} invariance~\cite[4.1]{TsiVerYor00} of~$\mcL^+_{\beta,\sigma}$ with respect to the action of the group of multipliers~$\exp(\mcC_c(X))\acts \Mbp(X)$ given by~$g.\eta\eqdef g\cdot \eta$ where~$g=e^h\in\mcC_b^+(X)$ for some~$h\in\mcC_c(X)$. The Radon--Nikod\'ym derivative satisfies in this case (see~\cite[4.1]{TsiVerYor00})
\begin{align*}
\forallae{\mcL^+_{\beta,\sigma}} \eta\in\Mbp(X) \qquad \frac{\diff (g.)_\pfwd \mcL^+_{\beta,\sigma}}{\diff \mcL^+_{\beta,\sigma}}(\eta)\equiv\exp (-\beta\,\sigma\ln g)\comm g\in \exp(\mcC_c(X)) \fstop
\end{align*}

\paragraph{The commutative action of~$\mfh_k$} It is the content of Theorem~\ref{c:Due}\iref{i:c:Due1} that the characteristic functionals of the measures~$\Dir{\boldalpha}$, varying~$\boldalpha\in\interior \Delta^{k-1}$, are projectively invariant under the action of the maximal toral subalgebra~$\mfh_k<\mfl_k$ in the representation~$\rho_\boldalpha$. Since~$\mfh_k$ acts on~$\mcO_\boldalpha$ by weight~$\boldalpha$ (see the proof of Thm.~\ref{c:Due}\iref{i:c:Due1}), for arbitrary~$J_\mbft\eqdef t_1J_{\alpha_1}+\cdots + t_kJ_{\alpha_k}\in \mfh_k$ one has
\begin{align*}
J_\mbft f_\boldalpha=(t_1\alpha_1+\cdots+t_k\alpha_k)f_\boldalpha=(\mbft \cdot \boldalpha) f_\boldalpha\qquad \mbft\in\R^k\fstop
\end{align*}

\paragraph{The non-commutative action of~$\mfl_k$ and a family of distinguished improper priors} In contrast to the case of the measures~$L_\boldalpha$ on affine spheres ---~where only the action of the commutative subgroup~$dSL^+_k(\R)<SL_k(\R)$ is taken into account~--- in the case of the Dirichlet distributions~$\Dir{\boldalpha}$ it is possible to detail the full non-commutative action of the algebra~$\mfl_k$ on their characteristic functionals.
Incidentally, let us notice that the acting object, although of special linear type in both cases, is a (subgroup of a) Lie \emph{group} in the first case, but the corresponding Lie \emph{algebra} in the latter case. This is because the action is, in the first case, an action on measures themselves, whereas, in the second case, on their characteristic functionals.

If~$\boldalpha\in\interior \Delta^{k-1}$, then~$(a)$ the action of basis elements in~$\mfr_k^+$ amounts to take (characteristic functionals of) Dirichlet-categorical \emph{posteriors}; it fixes the space~$\mcO_{H_{\boldalpha}}$ of (characteristic functionals of) such posteriors. On the other hand, $(b)$ the action of basis elements in~$\mfr_k^-$ amounts to take  (characteristic functionals of) Dirichlet-categorical \emph{priors}; such priors should be allowed to be improper, in the sense that they are no longer probability measures, but rather (in-)finite definite (i.e., positive or negative, \emph{not} signed) measures. Indeed, if we let~$\IDir{\boldalpha'}$ be any such improper prior, with density given by~\eqref{eq:DensityD} in the case when~$\boldalpha'\in \Lambda_\boldalpha^+$, then~$\IDir{\boldalpha'}$ has sign given by
\begin{align*}
\sgn (\Gamma(\boldalpha'))=\begin{cases}
1 & \text{if } \boldalpha'\in H_\boldalpha
\\
(-1)^{\ceiling{\alpha_1'}+\cdots + \ceiling{\alpha_k'}} & \text{otherwise}
\end{cases} \fstop
\end{align*}
The action of~$\mfr_k^-$ fixes the space~$\mcO_{\Lambda^+_\boldalpha}$ of (characteristic functionals of) all such priors and vanishes on the line~$M_{\boldalpha,0}$, the singular set of the normalization constant~$\Beta[\boldalpha']^{-1}$.
Finally, $(c)$ the action of basis elements in~$\mff_k$ contains every non-trivial combination of the actions~$(a)$ and~$(b)$, and fixes isoplethic hypersurfaces~$M_{\boldalpha,\ell}$, i.e. those where the intensity~$\boldalpha'$ has constant total mass~$\length{\boldalpha'}$.

In this framework, the case~$\boldalpha\in \bd\Delta^{k-1}$ is spurious, since the intensity measure~$\boldalpha$ should always be assumed fully supported.

\subsection{Infinite-dimensional statements}\label{ss:InfDimTDue}
For~$a\in \R$ we denote by~$\Mb^{>a}(X)$ the space of finite signed measures $\nu$ in~$\Mb(X)$ such that~$\nu X>a$.

\begin{thm}\label{t:Tre}
Let~$(X,\tau(X),\mcB(X))$ be a second countable locally compact Hausdorff space and~$\nu$ be a diffuse fully supported non-negative finite measure on~$X$. 
Let further\begin{align*}
\Phi[\nu,f]\eqdef& \sum_{n=0}^\infty \Poch{\nu X}{n}^{-1} Z_n(\nu f^1,\dotsc, \nu f^n)
\end{align*}
and
\begin{align*}
E_A\Phi[\nu,f]\eqdef& \int_A \diff\nu(y)\, \Phi[\nu+\delta_y,f] \comm\\
E_{A,-B} \Phi[\nu,f]\eqdef& \int_{A\setminus B} \diff\nu(y)\, \Phi[\nu+\delta_y,f] +\int_{B\setminus A}\diff\nu(y)\,\Phi[\nu-\delta_y,f]\fstop
%
\end{align*}

Then, 
\begin{enumerate}[$(i)$]
\item\label{i:t:Tre1} $\Phi[\nu,f]$ is a well-defined extension of the characteristic functional~$\FT{\DF_\nu}(f^*)$ on~$\Mb^{>0}(X) \times \mcC_c(X)$;
\item\label{i:t:Tre2} for every $\nu$ in~$\Mb^{>1}(X)$, every~$f$ in~$\mcC_c(X)$, every~$A,B$ in~$\mcB$, and every good approximation~$\seq{f_h}_h$ of~$f$ locally constant on~$\mbfX_h$ with values~$\mbfs_h$ for some~$\seq{\mbfX_h}_h\in\mfN\mfa(A,B\subset X)$, one has
\begin{align*}
E_A\Phi[\nu,f]=&\hlim \tonde{\sum_{i\mid X_{h,i}\subset A} E_{\alpha_{h,i}}} {}_{k_h}\!\Phi[\nu^\compo\mbfX_h,\mbfs_h]\comm\\
E_{A,-B}\Phi[\nu,f]=&\hlim \tonde{\sum_{\substack{i\mid X_{h,i}\subset A\setminus B\\ j\mid X_{h,j}\subset B\setminus A}} E_{\alpha_{h,i},-\alpha_{h,j}}} {}_{k_h}\!\Phi[\nu^\compo \mbfX_h,\mbfs_h] \comm
\end{align*}
where~$\boldalpha_h\eqdef \nu^\compo \mbfX_h$ and~$E_{\alpha_{h,i}}$, $E_{\alpha_{h,i},-\alpha_{h,j}}\in\mfg_{k_h}$.
\item\label{i:t:Tre3} let~$\sigma$ be a diffuse fully supported probability measure on~$(X,\tau(X))$ and let further~$\seq{\mbfX_h}_h\in \mfN\mfa(X,\tau(X),\sigma)$. For~$\sigma$-a.e.~$x$, such that~$X_{h,i_h}\downarrow_h \set{x}$, and for every good approximation~$\seq{f_h}_h$ of~$f$, locally constant on~$\mbfX_h$ and uniformly convergent to~$f$, there exist the pointwise limiting rescaled actions
\begin{align*}
\hlim \alpha_{h,i_h}^{-1}E_{\alpha_{i_h}} \FT{\DF_\sigma}(f_h^*)=& \FT{\DF_\sigma^x}(f^*)\comm\\
\hlim \alpha_{h,i_h}^{-1}J_{\alpha_{i_h}}=&\Id\comm\\
\hlim \alpha_{h,i_h}^{-1}E_{-\alpha_{i_h}}=&0\fstop
\end{align*}
\end{enumerate}

\begin{proof}
The functional~$\Phi[\nu, f]$ is well-defined in the first place since~$\nu X> 0$. For~\mbox{$c,t>0$} denote by~$P_{c,t}\subset \R^n$ the polydisk~$\set{\mbfy\in\R^n\mid \abs{y_i}\leq c t^i}$. By induction and~\eqref{eq:RecBellZ} it is not difficult to show that~$\max_{P_{c,t}}\abs{Z_n}=Z_n[c (t\uno)^{\compo \vec \mbfn}]$; moreover, by~\eqref{eq:MultiBell} and Theorem~\ref{t:MomDir}, the latter equals~$t^n\Poch{c}{n}/n!$. As a consequence, for arbitrary~$\nu$ in~$\Mb^{>0}(X)$ and~$f\in \mcC_c(X)$, letting~$y_i\eqdef \nu f^i$ above,
\begin{align*}
\abs{\Phi[\nu,f]}\leq& \sum_{n=0}^\infty \Poch{\nu X}{n}^{-1} \max_{P_{\norm{\nu}, \norm{f}}} \abs{Z_n}
=\sum_{n=0}^\infty \frac{\Poch{\norm{\nu}}{n}}{n!\Poch{\nu X}{n}} \norm{f}^n 
={}_1F_1\quadre{\norm{\nu}; \nu X; \norm{f}}\comm
\end{align*}
which is finite since~$\nu X>0$. This shows~\iref{i:t:Tre1}. Notably, if~$\nu$ is positive, then~$\abs{\Phi[\nu,f]}\leq \exp{\norm{f}}$ independently of~$\norm{\nu}$.

\smallskip

Let now~$A$ be in~$\mcB$ and $\seq{\mbfX_h}_h$ as in~\iref{i:t:Tre2}. Fix~$f$ in $\mcC_c(X)$, set~$\boldalpha_h\eqdef \nu^\compo \mbfX_h$ and let $\seq{f_h}_h$ be a good approximation of~$f$, locally constant on~$\mbfX_h$ with values~$\mbfs_h$.
Equation~\eqref{eq:t:dsa1} yields by summation
\begin{align}\label{eq:Borel}
\tonde{\sum_{i\mid X_{h,i}\subset A} E_{\alpha_i} } {}_{k_h}\!\Phi\quadre{\boldalpha_h;\mbfs_h}=\sum_{i\mid X_{h,i}\subset A} \alpha_{h,i} \,\, {}_{k_h}\!\Phi\quadre{\boldalpha_h+\mbfe_i;\mbfs_h} \fstop
\end{align}

More explicitly, since~$f_h$ is constant on each~$X_{h,i}$ with value~$s_{h,i}$, Proposition~\ref{p:LaurEGF} yields
\begin{align}\label{eq:t:Tre1}
\nonumber
\Bigg(&\sum_{i\mid X_{h,i}\subset A} E_{\alpha_i} \Bigg) {}_{k_h}\!\Phi\quadre{\boldalpha_h;\mbfs_h}=\\
\nonumber
&=\sum_{i\mid X_{h,i}\subset A} \nu X_{h,i} \sum_{n=0}^\infty \frac{1}{\Poch{\nu X+1}{n}} Z_n\tonde{\nu f_h+\frac{\nu(f_h\car_{X_{h,i}})}{\nu X_{h,i}}, \dotsc, \nu f_h^n+\frac{\nu(f_h^n \car_{X_{h,i}})}{\nu X_{h,i}}}\\
\nonumber
&=\sum_{n=0}^\infty \frac{1}{\Poch{\nu X+1}{n}} \times\\
\nonumber
&\quad\times\sum_{i\mid X_{h,i}\subset A} \int_{X_{h,i}}\!\!\!\!\diff\nu(y)\, Z_n\tonde{\nu f_h+\frac{\nu(f_h\car_{X_{h,i}})}{\nu X_{h,i}}\car_{X_{h,i}}(y), \dotsc, \nu f_h^n+\frac{\nu(f_h^n \car_{X_{h,i}})}{\nu X_{h,i}}\car_{X_{h,i}}(y) }\\
\nonumber
&=\sum_{n=0}^\infty \frac{1}{\Poch{\nu X+1}{n}} \sum_{i\mid X_{h,i}\subset A} \int_{X_{h,i}}\!\!\!\!\diff\nu(y)\, Z_n\tonde{\nu f_h+f_h(y), \dotsc, \nu f_h^n+f_h(y)^n }\\
\nonumber
&=\sum_{n=0}^\infty \frac{1}{\Poch{\nu X+1}{n}} \int_A\!\!\diff\nu(y)\, Z_n\tonde{\nu f_h+f_h(y), \dotsc, \nu f_h^n+f_h(y)^n } \fstop\\
&= \int_A\!\!\diff\nu(y) \sum_{n=0}^\infty \frac{1}{\Poch{(\nu+\delta_y) X}{n}} Z_n\tonde{\nu f_h+f_h(y), \dotsc, \nu f_h^n+f_h(y)^n } \fstop
\end{align}

Since~$\abs{f_h}\leq\abs{f}$ pointwise, the sequence~$\seq{f_h^i}_h$ converges strongly in~$L^1_\nu$ for every~$i\leq n$ for every~$n\in \N_1$, thus by continuity of~$Z_n$, there exists the limit
\begin{align*}
&\hlim \int_A\!\!\diff\nu(y) \sum_{n=0}^\infty \frac{1}{\Poch{\nu X+1}{n}} Z_n\tonde{\nu f_h+f_h(y), \dotsc, \nu f_h^n+f_h(y)^n }= E_A \Phi[\nu,f] \fstop
\end{align*}

The proof of the statement for~$E_{A,-B}$ is analogous. This completes the proof of~\iref{i:t:Tre2}. The requirement that~$\nu X>1$ is necessary to the convergence of~$\Phi[\nu-\delta_y,f]$ for~$y\in X$ in the definition of~$E_{A,-B}$, whereas it may be relaxed to~$\nu X>0$ in the case of~$E_A$. We will make use of this fact in the proof of~\iref{i:t:Tre3}.

\smallskip

Fix now~$x$ in~$X$ and let~$i_h\eqdef i_h(x)$ be such that~$X_{h,i_h}\downarrow_h \set{x}$. By~Lemma~\ref{l:MCNA}, the sequence~$\seq{i_h}_h$ is unique for~$\sigma$-a.e.~$x$. With the same notation of~\iref{i:t:Tre2}, let now~$A=X_{h,i_h}$ in~\eqref{eq:t:Tre1}. Then,
\begin{align}\label{eq:t:Tre2}
\alpha_{i_h}^{-1} E_{\alpha_{i_h}}\, {}_{k_h}\Phi[\boldalpha_h;\mbfs_h]=\alpha_{i_h}^{-1} E_{\alpha_{i_h}} \FT{\DF_\sigma}(f_h^*)=\frac{1}{\sigma X_{h,i_h}} \int_{X_{h,i_h}}\diff \sigma(y) \, \FT{\DF_{\sigma+\delta_y}}(f_h^*) \fstop
\end{align}

By~\eqref{eq:Thm1.2} and uniform convergence of the approximation
\begin{align}\label{eq:t:Tre3}
\hlim \abs{\tfrac{1}{\sigma X_{h,i_h}} E_{X_{h,i_h}}\Phi[\sigma,f_h]-\tfrac{1}{\sigma X_{h,i_h}} E_{X_{h,i_h}}\Phi[\sigma,f]}\leq \hlim e^{\norm{f}}\norm{f-f_h}=0\comm
\end{align}
thus,~\eqref{eq:t:Tre2} and~\eqref{eq:t:Tre3} yield, together with the continuity of~$y\mapsto \FT{\DF_{\sigma+\delta_y}}(f^*)$ for fixed~$f$ and~$\sigma$,
\begin{align*}
\hlim \alpha_{i_h}^{-1}E_{\alpha_{i_h}}\FT{\DF_\sigma}(f_h^*)=\hlim \frac{1}{\sigma_{X_{h,i_h}}}\int_{X_{h,i_h}}\diff\sigma(y)\, \FT{\DF_{\sigma+\delta_y}}(f^*)=\FT{\DF_{\sigma+\delta_x}}(f^*) \fstop
\end{align*}

By the Bayesian property~$\DF_\sigma^x=\DF_{\sigma+\delta_x}$, this yields the conclusion for the limiting raising action. 
Finally, since~$\sigma$ is a probability measure,~$\length{(\boldalpha_h)}=1$ for all~$h$, thus by Lemma~\ref{l:Commu},
\begin{align*}
\hlim \alpha_{i_h}^{-1} J_{\alpha_{i_h}} \FT{\DF_\sigma}(f_h^*)=&\hlim \FT{\DF_\sigma}(f_h^*)=\FT{\DF_\sigma}(f^*)\comm\\
\hlim \alpha_{i_h}^{-1} E_{-\alpha_{i_h}} \FT{\DF_\sigma}(f_h^*)=&\hlim 0=0\comm
\end{align*}
where the second equality for the first limiting action follows by~\eqref{eq:Thm1.1}.
In all three cases, independence of the limits from the chosen (good) approximation is straightforward. 
\end{proof}
\end{thm}

\section{Appendix}

We collect here some results in topology and measure theory.

\begin{lem}\label{l:MCNA}
Let~$(X,\tau(X),\mcB,\sigma)$ be a second countable locally compact Hausdorff Borel measure space of finite diffuse fully supported measure. Then, for every~$\seq{\mbfX_h}_h\in\mfN\mfa(X,\tau(X),\sigma)$ for $\sigma$-a.e.~$x$ in~$X$ there exists a unique sequence~$\seq{X_{h,i_h}}_h$, with~$i_h\eqdef i_h(x)$, such that~$\mbfX_h\ni X_{h,i_h} \downarrow_h\set{x}$.
\begin{proof}
Proposition~\ref{p:Alex} justifies well-posedness of the requirements in the definition of~$\seq{\mbfX_h}_h$.

Without loss of generality, each~$X_{h,i}$ may be chosen to be closed by replacing it with its closure~$\cl X_{h,i}=X_{h,i}\cup \bd X_{h,i}$. Hence~$\mbfX_h$ may be chosen to be consisting of closed sets (disjoint up to a $\sigma$-negligible set) with non-empty interior.
It follows by the finite intersection property that every decreasing sequence of sets~$\seq{X_{h,i_h}}_h$ such that~$X_{h,i_h}\in \mbfX_h$ admits a non-empty limit, which is a singleton because of the vanishing of diameters.
Vice versa, however chosen~$\seq{\mbfX_h}_h$, for every point~$x$ in~$X$ it is not difficult to construct a (possibly non-unique) sequence~$X_{h,i_h}$ (with~$i_h\eqdef i_h(x)$) convergent to~$x$ and such that~$X_{h,i_h}\in \mbfX_h$.
Furthermore, letting~$x$ be a point for which there exists more than one such sequence, we see that for every~$h$ the point~$x$ belongs to some intersection~$X_{h,i_{1}}\cap X_{h,i_{2}}\cap \dotsc$, hence, since every partition has disjoint interiors by construction,~$x\in \bd X_{h,i_1}\cap \bd X_{h,i_2}\cap \dotsc$.
Since for every~$h$ and~$i\leq k_h$ each set~$X_{h,i}$ is a continuity set for~$\sigma$, the whole union~$\cup_{h\geq 0}\cup_{i\in [k_h]} \bd X_{h,i}$ is $\sigma$-negligible, thus so is the set of points~$x$ considered above, so that for $\sigma$-a.e.~$x$ there exists a unique sequence~$\seq{X_{h,i_h}}_h$ such that~$X_{h,i_h}\in \mbfX^h$ and~$\hlim X_{h,i_h}=\set{x}$ and~$x$ belongs to each~$X_{h,i_h}$ in the sequence. 
\end{proof}
\end{lem}

Finally, recall the following form of L\'evy's Continuity Theorem.
\begin{thm}[{\cite[Thm.~3.1, p.~224]{VakTarCho87}}] Let~$(Y,\tau(Y))$ be a completely regular Hausdorff topological space,~$V$ be a linear subspace of~$\mcC(Y)$ separating points in~$Y$ and~$\chi$ be a complex-valued functional on~$V$. If~$\seq{\mu_\gamma}_\gamma$ is a narrowly precompact net of Radon probability measures on~$(Y,\mcB(Y))$ and~$\lim_\gamma \FT{\mu_\gamma}(v)=\chi(v)$ for every~$v$ in~$V$, then~$\seq{\mu_\gamma}_\gamma$ converges narrowly to a Radon probability measure~$\mu$, the characteristic functional thereof coincides with~$\chi$.
\end{thm}

\begin{cor}\label{c:Vak} Let~$\seq{\mu_\gamma}_\gamma$ be a narrowly precompact net of random probabilities over the space~$(X,\mcB(X))$. If $\lim_\gamma \FT{\mu_\gamma}(f^*)=\chi(f^*)$ for every~$f$ in~$\mcC_c(X)$, then~$\seq{\mu_\gamma}_\gamma$ converges narrowly to a random probability~$\mu$, the characteristic functional thereof coincides with~$\chi$.
\begin{proof}
By Proposition~\ref{p:Alex} the space~$(X,\tau(X))$ is Polish, hence so is~$\Mbp(X)$~\cite[15.7.7]{Kal83}, thus the space $\Meas^+_{\leq 1}(X)\eqdef \set{\mu\in \Mbp(X)\mid \mu X\leq 1}$ is too, being closed, and~$\msP(X)$, being a $G_\delta$-set in~$\Meas^+_{\leq 1}(X)$. Since every finite measure on a Polish space is Radon~\cite[Thm.~7.1.7]{Bog07}, each~$\mu_\gamma$ is Radon. Consider~$\Mb(X)$ endowed with the vague topology. The dense subset~$\mcC_c(X)$ of the topological dual $(\Mb(X),\tau_v(\Mb(X)))'=\mcC_0(X)$ separates points in~$\Mb(X)$, hence it separates points in~$\msP(X)\subset \Mb(X)$. The conclusion follows now by the Theorem choosing~$Y=\msP(X)$ and~$V=\mcC_c(X)$.
\end{proof}
\end{cor}



\end{document}